\documentclass{article}

\usepackage{graphicx}
\usepackage{amssymb,amsmath}
\usepackage[english]{babel}
\usepackage{color}
\newcommand{\Gn}{\ensuremath{\boldsymbol \Gamma_{\!\rm{n}}}}
\newcommand{\Gpr}{\ensuremath{\boldsymbol \Gamma_{\!\rm{pr}}}}
\newcommand{\Gpo}{\ensuremath{\boldsymbol \Gamma_{\!\rm{pt}}}}
\newcommand{\dn}{\ensuremath{\mathbf d}}
\newcommand{\nnu}{\ensuremath{\boldsymbol \nu}}
\newcommand{\nnumap}{\ensuremath{\nnu_{\!\rm{\scriptscriptstyle MAP}}}}

\begin{document}

\title{Quantifying uncertainty in inverse scattering problems set
in layered environments 
}

\author{Carolina Abugattas, Ana Carpio, Gerardo Oleaga, \\
Universidad Complutense de Madrid, Spain \\[1ex]
Elena Cebrián, Universidad de Burgos, Spain}

\maketitle

{\bf Abstract.} 
The attempt to solve inverse scattering problems often leads to optimization 
and sampling problems that require handling moderate to large amounts of 
partial differential equations acting as constraints.
We focus here on determining inclusions in a layered medium from the 
measurement of wave fields on the surface, while quantifying uncertainty 
and addressing the effect of wave solver quality.
Inclusions are characterized by a few parameters describing their material 
properties and shapes. We devise algorithms to estimate the most likely 
configurations by optimizing cost functionals with Bayesian regularizations
and wave constraints.
In particular, we design an automatic Levenberg-Marquardt-Fletcher type 
scheme based on the use of algorithmic differentiation and adaptive finite 
element meshes for time dependent wave equation constraints with changing 
inclusions. In synthetic tests with a single frequency, this scheme converges 
in few iterations for increasing noise levels.
To attain a global view of other possible high probability configurations 
and asymmetry effects we resort to parallelizable affine invariant Markov 
Chain Monte Carlo methods, at the cost of solving a few million wave problems.
This forces the use of prefixed meshes. While the optimal configurations
remain similar, we encounter additional high probability inclusions influenced 
by the prior information, the noise level and the layered structure, effect that
can be reduced by considering more frequencies.



\section{Introduction}
\label{sec:intro}

Geophysical imaging is a noninvasive imaging technique that investigates 
the subsurface. It uses elastic waves created by man-made explosions and 
vibrations  to image the subsurface of Earth \cite{seismic_imaging}.
Elastic waves split in two components \cite{landau,seismic_waves}.
Compression waves (P-waves) are longitudinal in nature and penetrate 
through subsurface layers, causing the ground to compress and stretch along 
the axis of propagation of the wave in a similar way to sound waves. 
Due to these properties, longitudinal P-waves constitute a useful tool 
for underground imaging at considerable depth.

\begin{figure}[!hbt]
\centering
\includegraphics[width=8cm,angle=0]{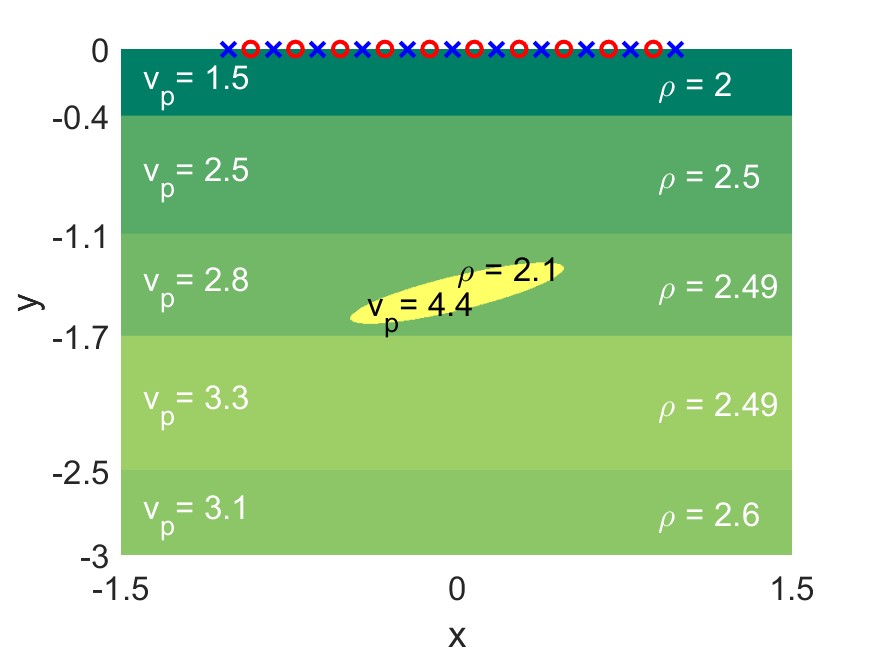} 
\caption{Schematic representation of an imaging set-up. The emitters (red)
generate waves which interact with the medium. The reflected waves
are recorded at the receivers (blue). The velocities are in units of $10^3$
{\tt m/s} and the densities in units of {\tt kg/m}$^3$. Units for $x$ and $y$ 
are {\tt km}. Parameter values are typical of sandstone, shale, limestone,
and salt for the inclusion.
}
\label{fig1}
\end{figure}

Figure \ref{fig1} represents the structure of a standard imaging set-up. 
Explosions or strong impacts at certain locations produce elastic waves that 
propagate under the surface. These waves interact with subsurface structures 
and are reflected at different depths 
\cite{fichtner_travel,seismic_waves,seismic_imaging,papanico}.  
Then, a grid of recording devices measures the reflected waves. Different 
procedures have been developed to image subsurface layer properties, such 
as classical seismic reflection and full waveform inversion 
\cite{seismic_imaging, virieux}.
Assuming that the basic layered structure has been already characterized by 
other techniques \cite{seismic_waves,georg_migration}, our goal here is to 
refine the description of localized inclusions of different materials with quantified 
uncertainty, which leads to solving an inverse scattering problem.

Inverse scattering problems share a similar mathematical structure: a number of
sources launch waves that interact with a medium containing scatterers and the 
resulting waves are recorded at a set of detectors. Knowing the recorded signals 
we aim to identify the scatterers. 
Methods are adapted to the nature of the waves and the way the inclusion 
properties are represented mathematically, see
\cite{park,greengard,cakonilsm,jcp19,conca,dorn,feijoo,laurain,
hohage2d,kirschmusic,litman} and references therein, for example.
Waves can be time-harmonic, and thus governed by stationary elliptic problems, 
or time dependent (thermal, electromagnetic, acoustic...).
Physical properties can be represented as infinite dimensional coefficient functions, 
high dimensional sets of values at grid points, level set functions or localized 
star-shaped inclusions, for instance.
Inverse scattering problems are severely ill-posed and are often regularized
by means of variational formulations \cite{colton,feijoo}.
A variety of deterministic approaches are based on minimizing cost functionals 
which compare the recorded data with the synthetic data that would be obtained
for arbitrary scatterers, as predicted by a selected forward model for 
the propagation of the emitted waves. 
Optimization strategies employed often implement total variation regularizations, 
Tikhonov regularizations or iteratively regularized Gauss Newton approaches
\cite{jcp19,chan,total1,blocky,hohage2},  and may rely on distances other than 
the Euclidean, i.e. Wasserstein distances \cite{yunan}.

In practice, recorded data are affected by noise. Typically, deterministic strategies 
produce the best solution for a given dataset, that is, a given realization of
noise. However, we lack information on how the proposed solution can change
as we vary the noise realizations. Under some conditions, deterministic iteratively 
regularized Gauss Newton methods can provide approximate solutions that are 
acceptable within a noise level \cite{hohage2}. 
Instead, Bayesian formulations of the inverse problem aim to provide the 
most likely solutions while quantifying uncertainty about them for a given noise 
magnitude. Statistical inverse problems also require the choice of a mathematical 
representation of the unknown fields, see
\cite{buishape,elena,matt_thesis,capistran,georg_migration} for high dimensional 
sets of values at grid points, level set functions or localized star-shaped inclusions,
for instance. In a Bayesian context, the problem can be addressed by Markov
Chain Monte Carlo (MCMC) techniques or Laplace approximations 
\cite{mcmc,inference}.  In contrast with  infinite dimensional formulations
\cite{buishape,georg_linearized,matt,stuart}, we will focus here on situations in 
which the inclusions are characterized by a finite and fixed number of random 
variables. Our goal is to devise automatic procedures and assess specific issues 
related to the practical choice of time dependent wave equation solvers. Furthermore, 
this choice allows us to study at a reasonable computational cost the appearance 
of several high probability configurations and how we can modify the 
imaging set up  to eliminate spurious ones, reducing uncertainty in this way.

Success of both deterministic and statistical approaches to the numerical solution 
of inverse scattering problems relies on being able to solve moderate to large 
amounts of auxiliary boundary value problems modeling wave propagation. 
The configuration under here study (see Figure \ref{fig1}) shows abrupt changes at 
interfaces between materials of a different nature, which are described as discontinuities 
in the fields representing physical properties. In time-harmonic problems, efficient
boundary element solvers are available \cite{jcp19,harbrech,hohage2d}, which
are usually applied to inclusions in homogeneous backgrounds.
For time dependent problems, boundary element methods are less developed
\cite{tonatiuh}. Moreover, in our layered geometry, interfaces arise not only
around the sought inclusions, but in the surrounding environment too.
When we employ finite element  (FEM) solvers for wave equations of the form
$\rho(\mathbf x) u_{tt} -{\rm div}(\chi(\mathbf x) \nabla u) = h$, 
with piecewise fields $\rho$ and $\chi$ representing the density and
elastic constants of different materials, the question arises of whether to keep 
a fixed mesh (which increases the error at discontinuities) or adapt the mesh 
to the varying interfaces (which increases the computational  cost) and whether 
this choice has a relevant effect on the results in our context, since it may affect 
the speed of propagation when crossing layers 
\footnote{Some work in the literature, such as \cite{matt}, use wave constraints
of the form  $m(\mathbf x) u_{tt} - \Delta  u = h$ thinking 
of $m^{-1/2}=v_{\rm p}$ as the wave velocity, instead of 
$\rho(\mathbf x) u_{tt} -{\rm div}(\chi(\mathbf x) \nabla u) = h$.
When the density $\rho$ and elastic parameters $\chi$ are constant everywhere,
we can indeed write such an equation for the waves, with constant $m$ and
$v_{\rm p}= \sqrt{\chi /\rho}$, see Appendix.  If $\chi$ varies spatially, we cannot.
Furthermore, in medical applications, for instance, we can set $\rho \sim 1$ and
work with one field instead of two \cite{elena}. In geophysical applications,
we must keep two fields, either $\rho(\mathbf x)$ and $\chi(\mathbf x)$ or
$\rho(\mathbf x)$ and $\rho(\mathbf x) v_{\rm p}(\mathbf x)^2$ \cite{georg_linearized}.
Reference \cite{matt} analyzes the performance of different distances
in the costs/likelihoods in uncertainty studies in a two dimensional geometry
similar to ours, but working in an infinito-dimensional framework with a more or
less homogeneous background. However, their choice of wave constraint
$m u_{tt} - \Delta  u = h$  removes one of the two fields characterizing the media
and the inclusions.}.

In this paper, we develop techniques to estimate the most likely inclusions 
given noisy data, characterized as minima of cost functionals constrained 
by time dependent wave equations with a Bayesian regularization. 
Regularizing terms are expected to convexify the cost functionals, still,
local minima may persist depending on the quality of the prior information 
and the imaging set-up design.
We consider three approaches, which provide complementary information.
First, we propose an automatic adaptive optimization scheme that combines 
Levenberg-Marquardt  \cite{lmf, lm} type iterations with the use of adaptive 
meshes \cite{Persson,Strang-Persson} for finite element discretization  
of the wave constraint and  automatic differentiation  \cite{ForwardDiff} to 
quantify variations in the proposed inclusions at each stage. 
The resulting adaptive  scheme is fully automatic: the user must only provide 
the data, the noise level, and the prior information. 
This algorithm usually converges to the main minimum in a few iterations
in the tests we have performed.
Second, we introduce an alternative scheme based on FEM discretizations
of the wave constraint on fine enough prefixed stratified or uniform meshes. 
This algorithm employs finite difference approximations of the observation 
operator defined by the wave problems as we modify the inclusions to seek 
descent directions. Varying the step, we reach the global minimum but could
also identify some additional local minima if present.
Finally, we resort to affine invariant MCMC samplers \cite{goodman} to 
characterize all minima in detail. Adaptive meshes are currently not affordable
within MCMC schemes, therefore we employ again fixed uniform or stratified
meshes. The optimal configurations encountered by the three methods are
similar in the synthetic tests we have carried out. 
We obtain basic uncertainty estimates on parameter ranges resorting to the 
so-called Laplace approximation, by linearizing about the optimal inclusions.
Unlike what happens in similar imaging set-ups for medical applications in 
homegeneous backgrounds \cite{elena} when employing one frequency, 
here MCMC techniques identify several  high probability configurations, regardless 
of the use of uniform (blind to the layered structure) or stratified (adapted to the 
known layered structure but not to the scatterer) meshes. 
Provided the meshes are fine enough, 
their effect on the hight probability configurations remains small. The highest 
probability estimate resembles the true inclusion. Secondary configurations seem 
to reflect the interaction with the layered structure, a memory of the prior information, 
the difficulty to resolve depth in this type of imaging set-ups with one sided 
information and the noise to signal ratio. The appearance of `phantom' inclusions in 
uncertainty studies has been related to aberrations caused by lack of information 
from enough incidence directions \cite{sergei} in optical applications. Here, we notice
that the use of data recorded for different frequencies helps to convexify the cost 
and reduces the uncertainty caused by secondary modes.

The paper is organized as follows. We describe the imaging set-up and state
the inverse problem in Section \ref{sec:true_ip}.  
Section \ref{sec:approximate_ip} formulates finite dimensional approximations
obtained discretizing the forward problem for wave propagation and the observation 
operator. In Section \ref{sec:optimization}, we develop optimization approaches to 
estimate the most likely inclusions given noisy data.
After testing the performance of the algorithms with synthetic data, in Section 
\ref{sec:uncertainty} we quantify uncertainty in the inclusion parameter predictions 
by Laplace approximations and by MCMC sample analysis. We choose a practical
test of current interest in geological studies: the characterization of a salt inclusion,
which usually hide reservoirs of raw materials \cite{dome}.
Lastly, Section \ref{sec:conclusions} summarizes our conclusions.  Appendix A collects 
details on the nondimensionalization of the model and the parameters selected.
Appendix B formulates a well posedness theory for the truncated wave problem 
with artificial boundary conditions we use as constraint. Finally, Appendix C 
describes the numerical schemes we implement to  approximate its solutions and 
discusses their convergence properties.


\section{Object based full-waveform inversion}
\label{sec:true_ip}

Consider an imaging set-up as depicted in Figure \ref{fig1}. A grid of sources,
$\mathbf x_k$, $k=0,\ldots,K$, emits waves that propagate under the ground, 
interact with subsurface structures, are reflected and recorded again on the
surface, at an interspaced grid of receivers $\mathbf r_j$, $j=0,\ldots,J$.
Here, we assume that we have information on the layered structure. The goal
is to characterize localized inclusions of different materials given the data 
$\mathbf d$ recorded at the receivers at a sequence of times $t_m$, 
$m=1,\ldots,M$. These scatterers are often described by a finite  collection of 
parameters $\nnu$ representing their material parameters and geometry 
\cite{babak,jcp19,blocky,harbrech,capistran}. To identify $\nnu$
from measurements $\mathbf d$ we first need to relate both by means of a 
mathematical model of wave propagation and a suitable observation operator.


\subsection{Observation operator}
\label{sec:obs_op}

In a simple model of longitudinal wave propagation, the waves emitted by the
sources are governed by a scalar wave equation:
\begin{eqnarray} \begin{array}{ll}
\rho \, u_{tt} = {\rm div}[(\lambda +2 \mu)
\nabla u] +\rho(\mathbf x)  h & \mathbf x \in {\cal R}, \, t > 0, \\ [1ex]
\nabla u \cdot \mathbf n = 0,  & \mathbf x \in \Sigma, \, t > 0, \\ [1ex]
u(0, \mathbf x) = 0, \; u_t(0, \mathbf x) =0,  & \mathbf x \in R.
\end{array} \label{balance_half} 
\end{eqnarray}
Here, the subsurface is represented by a half-space ${\cal R}$ and
$u(\mathbf x,t)$ denotes downward displacements in the direction 
$y$, with $\mathbf x = (x,y)$.
We model the emitters as source terms of the form $f(t) g(\mathbf x - 
\mathbf x_k)$, $k=1,\ldots, K$. The function $f(t)$ is a Ricker wavelet $f(t)=f_0 
(1- 2 \pi^2  f_M^2 t^2) e^{-\pi^2 f_M^2 t^2}$ with peak frequency $f_M$. The 
function $g$ has a narrow support  and zero normal derivative at $y=0$. For 
instance, we may set  $h(t, \mathbf x) = f(t) G(\mathbf x ) =  \frac{f_0 f(t)}
{( \pi \kappa )^{n/2} } \sum_{k=1}^K \exp(- \frac{|\mathbf x - \mathbf x_k|^2}
{\kappa })$, $n=2$.

In a layered geometry, the density $\rho(\mathbf x)$ and the elastic constants 
$\lambda(\mathbf x)$ and $\mu(\mathbf x)$ are piecewise constant in ${\cal R}$.
When additional inclusions $\Omega_{\nnu} = \cup_{\ell=1}^L  \Omega^\ell$ of 
different materials are present we set
\begin{eqnarray}  
\rho(\mathbf x) = \left\{ \begin{array}{ll} 
\rho,    & \mathbf x \in {\cal R} \setminus \overline{\Omega}_{\nnu}, \\
\rho_{\rm i}^\ell, & \mathbf x \in  \Omega^\ell, \; \ell = 1, \ldots, L,
\end{array}  \right.  \label{coefrho}  \\
\chi(\mathbf x) = \left\{ \begin{array}{ll} 
\lambda + 2 \mu  = \rho v_{\rm p}^2,    & \mathbf x \in {\cal R} \setminus
\overline{\Omega}_{\nnu}, \\
\lambda_{\rm i}^\ell + 2 \mu_{\rm i}^\ell = \rho_{\rm i}^\ell (v_{\rm p, \rm i}^2)^\ell,
&  \mathbf x \in  \Omega^\ell,  \ell = 1, \ldots, L,
\end{array} \right.
\label{coefs}
\end{eqnarray}
where $v_{\rm p}$ is the wave speed in each of the subdomains.

Measurements can be related to wave models in different ways depending on the
set-up. Here we assume that the measured data correspond to values taken
by the solution of (\ref{balance_half}) at the receivers. Given inclusions characterized
by parameters $\nnu_{\rm true}$, the corresponding data $\mathbf d_{\rm true}$ are 
related to $\nnu_{\rm true}$ through the observation operator
\begin{eqnarray} \begin{array}{ll}
\mathbf o: & \mathbb R^P \longrightarrow \mathbb R^D 
\\ [1ex]  
&\nnu \longrightarrow (u_{\nnu}(\mathbf r_j, t_m))_{j=1,\ldots,J, m=1,\ldots,M},
\end{array} \label{observation} 
\end{eqnarray}
where $u_{\nnu}$ stands for the solution of (\ref{balance_half}), $P$ is the number
of parameters and the number of data is $D=MJ$.
We expect the data $\mathbf d_{\rm true}$ recorded for inclusions parametrized 
by $\nnu_{\rm true}$ to satisfy $\mathbf o(\nnu_{\rm true}) = \mathbf d_{\rm true}$.
In practice, the recorded data are corrupted by noise:
\begin{eqnarray}
\mathbf d = \mathbf o(\nnu_{\rm true}) + \boldsymbol \varepsilon, \quad
\boldsymbol \eta \sim {\cal N}(0,\Gn).
\label{druido}
\end{eqnarray}
We assume that $\boldsymbol \varepsilon$ is distributed as a multivariate Gaussian.


\subsection{Deterministic inverse problem}
\label{sec:det_ip}


Given data $\dn$, the inverse problem consists in finding the parameters $\nnu$
for which the observed data $\dn=\mathbf o(\nnu)$. This is a severely ill posed 
problem  \cite{colton} which is often regularized by means of constrained optimization
reformulations aiming to optimize a cost \cite{pattern} 
\begin{eqnarray}
J(\nnu) = J_{\rm d}(\nnu)+ {\cal R}(\nnu), \quad 
J_{\rm d}(\nnu) = {1\over 2} \sum_{j=1}^J \sum_{m=1}^M 
|u_{\nnu}(r_j,0,t_m) - d_j^m|^2,
\label{dcost}
\end{eqnarray}
where $u_{\nnu}$ is the solution of the forward problem (\ref{balance_half}) for 
inclusions defined by parameters $\nnu$. To this purpose, iterative optimization 
techniques are often employed \cite{lmf,lm,nocedal}. In the absence of noise,
the parameters characterizing the exact inclusions furnish a global minimum of 
the deterministic cost $J_d(\nnu)$. Additional local minima may be present 
depending on the data quality. The term $R(\nnu)$ is a regularizing term aiming
to convexify the cost and ensure the occurrence of a unique global minimum. 
Typical choices are Tikhonov and total variation terms \cite{chan,total1} but also 
iteratively regularizing Gauss Newton schemes, in which the regularizing term 
vanishes while optimizing \cite{hohage2}. The later schemes have the potential 
of providing solutions which are robust for certain noise levels.


\subsection{Bayesian inverse problem}
\label{sec:bay_ip}

Bayesian approaches aim to quantify uncertainty in the solution of the inverse
problem relying on Bayes' formula \cite{KaipioSomersalo06,Tarantola05}. In finite 
dimension,  we consider the parameters random as variables with posterior
density
\begin{eqnarray}
p_{\rm pt}(\nnu):= p(\nnu|\dn) = \frac{ p(\dn|\nnu)} {p(\dn)} p_{\rm pr}(\nnu), 
\label{bayes}
\end{eqnarray}
given data $\dn$. The prior density $p_{\rm pr}(\nnu)$ represents the available 
prior knowledge whereas the conditional probability $p(\dn|\nnu)$ describes the 
likelihood of the measurements $\dn$ given the parameters $\nnu$. The 
density $p(\dn)$ is a normalization factor, to keep the integral of the posterior
probability equal to one. Here, we choose a likelihood
\begin{eqnarray}
p( \dn | \nnu) = \frac{1}{(2\pi)^{N/2} \sqrt{|\Gn|}} \exp \left(- \frac{1}{2} \|
\mathbf o( \nnu) - \dn  \|^2_{\Gn^{-1}} \right),
\label{likelihood} 
\end{eqnarray}
where $\| \mathbf v \|_{\Gn^{-1}}^2 =  \mathbf {\overline v}^t \Gn^{-1} \mathbf v$,
with $t$ denoting transpose,
$D$ is the dimension of $\mathbf d$  and $\mathbf o(\nnu)$ represents the 
observation operator. 
We consider a diagonal covariance matrix $\Gn$ with constant diagonal of 
magnitude $\sigma_{\rm noise}^2$. A typical choice for $p_{\rm pr}(\nnu)$ is a 
multivariate Gaussian with covariance matrix $\Gpr$ 
\begin{eqnarray} 
p_{\rm pr}(\nnu)  = \left\{ \begin{array}{ll}
\frac{1}{(2 \pi)^{P/2}}  \frac{1}{\sqrt{|\Gpr |}} \exp \left(-\frac{1}{2}
\|\nnu - \nnu_0\|_{\Gpr^{-1}}^2 \right) & \nnu \in {\cal P}, \\
0 & \nnu \notin {\cal P},
\end{array} \right.  \label{prior}
\end{eqnarray}
where ${\cal P}$ represents a set of constraints to be satisfied by $\nnu$.
The resulting posterior probability is proportional to 
\begin{eqnarray}
p(\nnu |  \dn ) \sim \exp \left( -\frac{1}{2} \| \mathbf o( \nnu) - \dn  \|^2_{\Gn^{-1}}
-\frac{1}{2} \|\nnu - \nnu_0\|_{\Gpr^{-1}}^2 \right), \quad \nnu \in {\cal P}.
\label{unnormalized}
\end{eqnarray} 
The solution of the inverse Bayesian problem is the characterization of
the posterior probability. Full characterization of this unnormalized posterior 
distribution is a challenging problem that can be addressed by Markov Chain 
Monte Carlo (MCMC) methods \cite{pattern,mcmc} to a certain extent.
The most likely set of parameter values $\nnumap$ defines the MAP (maximum 
a posteriori) estimate and minimizes the functional
\begin{eqnarray}
J(\nnu) = \frac{1}{2} \| \mathbf o( \nnu) - \dn  \|^2_{\Gn^{-1}}
+ \frac{1}{2} \|\nnu - \nnu_0\|_{\Gpr^{-1}}^2,
\label{bcost}
\end{eqnarray}
which is  expected to be convex  for good choices of  $\nnu_0$ 
and $\Gpr$, however, the prior information might not be good enough.
Compared to (\ref{dcost}), the first term is the deterministic
cost $J_{\rm d}(\nnu)$ scaled by $\sigma_{\rm noise}^2$, while the second
acts as  a Bayesian regularizing term ${\cal R}(\nnu)$.
Once we have calculated $\nnumap$, the so-called Laplace approximation
linearizes the posterior distribution about $\nnumap$ and estimates the
posterior  distribution through a multivariate Gaussian 
distribution  ${\cal N}(\nnumap, \Gpo)$ \cite{inference,georg_mcmc}.
This may be a useful approximation when the posterior distribution is not
multimodal, otherwise it only captures the main mode.


\section{Approximate inverse problem}
\label{sec:approximate_ip}

For computational studies, the continuous observation operator (\ref{observation})
must be replaced by a discrete approximation. The accuracy of this approximation
depends on the scheme used to construct numerical solutions of (\ref{balance_half}).
In homogeneous backgrounds, and for stationary constraints, boundary value elements
furnish an efficient procedure that minimizes numerical artifacts \cite{jcp19,harbrech,hohage2}
and computational costs associated to mesh design.
In a layered medium and with time dependent constraints, we are forced to rely on
finite element methods and face issues related to mesh design and computational cost.


\subsection{Truncated forward problem}
\label{sec:balance_truncated}

Due to the finite speed of propagation of waves, the wave field solution of (\ref{balance_half})
will only be nonzero inside the domain of influence of the emitted waves  \cite{john}. 
One can exploit this fact to truncate the computational domain to a finite rectangular
domain. Given a time $\tau>0$, we can truncate the half-space to a large box $R_\tau$
in such a way that the wave field vanishes at its bottom and lateral boundaries during
the time interval $[0,\tau]$. Then, the problem set in the whole halfspace is 
equivalent to the problem set in $R_\tau$ with zero Dirichlet or Neumann boundary 
conditions at those boundaries for $t \in [0,\tau]$. 
For computational purposes, we often need the computational region to be
as small as possible. Under some conditions the original model can be replaced by
an equivalent problem set in a smaller rectangle $R$ provided adequate nonreflecting 
boundary conditions are available. A typical choice for scalar wave equations with 
constant coefficients are  conditions of the form ${\partial u \over \partial \mathbf n} 
= - {u_t \over v_{\rm p}}$, where $v_{\rm p}$ is the wave speed, see \cite{nonreflecting}. 
Similar conditions can be exploited in layered geometries, as depicted in Figure \ref{fig1}.
The truncated problem is then
\begin{eqnarray} \begin{array}{ll}
\rho(\mathbf x) \, u_{tt} = {\rm div}(\chi(\mathbf x) \nabla u) +
\rho(\mathbf x)  h(t,\mathbf x), & \mathbf x \in R, \, t \in [0,T], \\ [1ex]
\nabla u \cdot \mathbf n = 0,  & \mathbf x \in \Sigma, \, t \in [0,T], \\ [1ex]
\nabla u \cdot \mathbf n = - \gamma(\mathbf x) u_t,  & \mathbf x \in 
\partial R \setminus \overline{\Sigma},  \, t \in [0,T], \\ [1ex]
u(0, \mathbf x) = u_0(\mathbf x), \; u_t(0, \mathbf x) = u_1(\mathbf x),  & 
\mathbf x \in R,
\end{array} \label{balance_escalar} 
\end{eqnarray}
where $\rho$ and $\chi$ are defined in (\ref{coefrho}) and (\ref{coefs})
and $\gamma= \rho v_{\rm p}$. 
Existence of solutions to this problem when $\gamma >0$ is not  immediate, since the 
boundary condition uses values of $u_t$ on the boundary. Standard existence results 
for wave equations  guarantee the existence of solutions with $L^2(R)$ derivatives $u_t$. 
However, for $u_t$ to have a trace on $\partial R$  we would need at least $H^1(R)$ 
regularity. Appendix B establishes existence, uniqueness, regularity and stability results 
for (\ref{balance_escalar}). When $\gamma=0$, we have a standard Neumann problem.


\subsection{Choice of mesh and discretization}
\label{sec:balance_mesh}

\begin{figure}[!htb] \centering
\includegraphics[width=12cm]{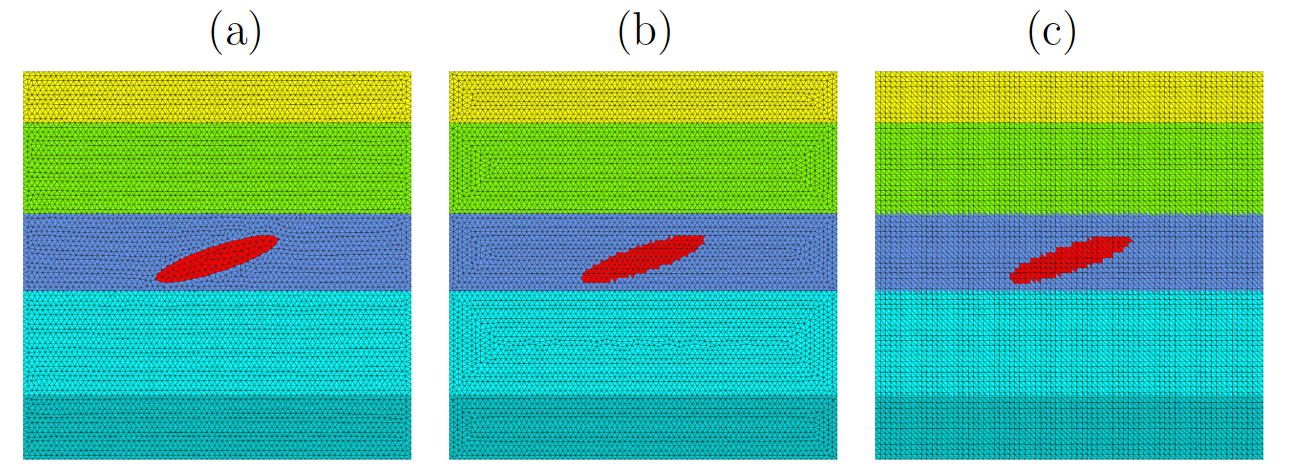} 
\caption{Types of meshes considered: 
(a) adapted to all the  subdomains, 
(b) adapted to the stratified structure but not to the changing inclusion, 
(c) uniform.}
\label{fig2}
\end{figure}

Problem (\ref{balance_escalar}) is set in layered domains with inclusions represented
by piecewise constant coefficients. Domain decomposition techniques \cite{bamberger}
provide effective tools to address it. For this purpose, the spatial mesh must adapt to 
the interfaces: all the triangles must be entirely contained in one subdomain and
connecting triangles share vertices lying on the boundaries between domains, see
Figure \ref{fig2}(a).
Studying the posterior distribution (\ref{unnormalized}) by MCMC techniques requires
solving millions of forward problems with different subdomains, as defined
by different $\nnu$ choices. Domain decomposition techniques require new meshes
for each choice of $\nnu$ and become unaffordable, forcing the use of rougher 
approximations built on fixed meshes, see Figure \ref{fig2} (b) and (c).
In contrast, optimization problems (\ref{dcost}) and (\ref{bcost}) can be addressed with 
either adaptive or fixed meshes, which allows us to evaluate the effect of using different 
meshes, while keeping the same discretization schemes.

For a given mesh, Appendix C describes the discretization procedure we have employed 
in the tests performed here, as well as their approximation properties.
Figure \ref{fig2} compares the effect of the triangulation choices on the numerically 
observed data at the receivers. The largest errors are observed when comparing results 
obtained with the adapted and uniform meshes, while the smallest errors correspond
to comparison of results calculated with the adapted and stratified mesh (adapted to
the layers, but not to the inclusion). As the maximum diameter of the triangulation 
elements $\delta x$ tends to zero, we  can prove convergence for adapted meshes,
see Appendix C. A better approximation of the transmission conditions for wave 
propagation at the interface between subdomains is expected.

\begin{figure}[!htb]
\centering
\hskip -4.2cm (a) Data \hskip 5.5cm (b) Error \\
\includegraphics[trim= 2.1cm 5mm 2mm 5mm, clip=true,width=6.0cm]{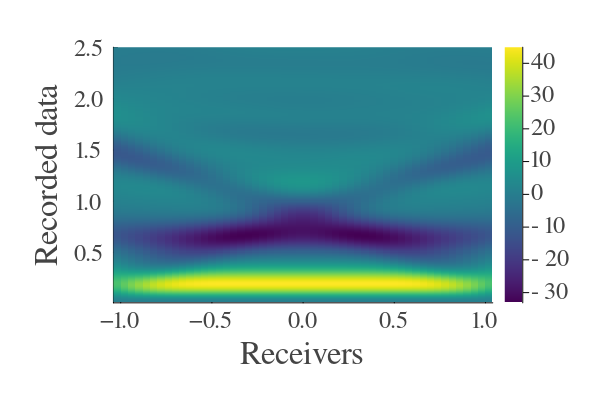}
\includegraphics[trim= 2.1cm 5mm 2mm 5mm, clip=true,width=6.0cm]{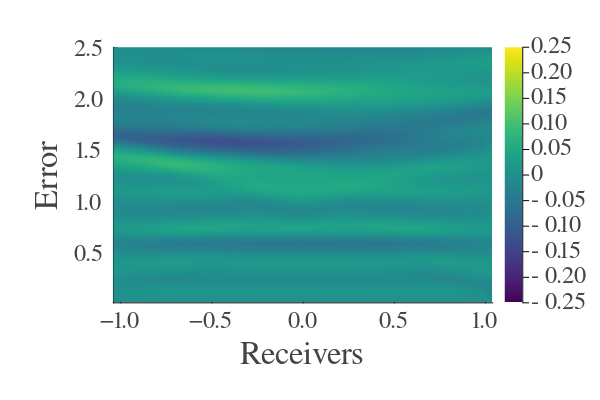} \\
\hskip -4.2cm (c)  Error \hskip 5.5cm (d)  Error \\
\includegraphics[trim= 2.1cm 5mm 2mm 5mm, clip=true,width=6.0cm]{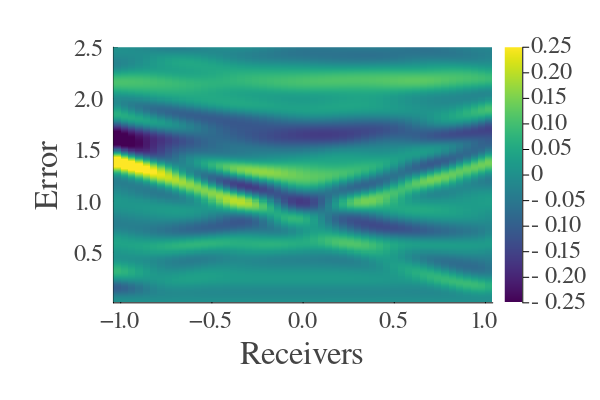}%
\includegraphics[trim= 2.1cm 5mm 2mm 5mm, clip=true,width=6.0cm]{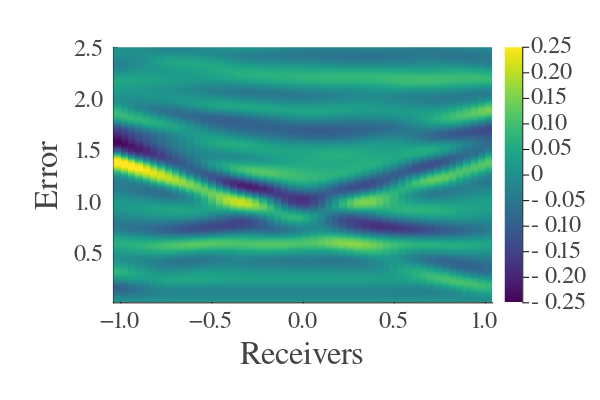}%
\caption{
(a) Values recorded at the receivers with the adaptive mesh in Figure \ref{fig2}(a)
for the parameter values  and geometry  in Figure \ref{fig1} with 
$\delta x = 0.04 $, $\delta t = 1\text{e-}3$. 
Profiles are recorded at intervals of $0.01$.
Errors when comparing the data in (a) with  the observed values obtained with 
(b) the stratified mesh in Figure \ref{fig2}(b) and 
(c) the uniform mesh in Figure \ref{fig2}(c) and
when comparing the data obtained with the latter two meshes between themselves 
(d).
}
\label{fig3}
\end{figure}

The numerical tests presented in this paper use synthetic data generated employed
either of the three meshes and corrupted by noise up to a 15\% level to solve
a nondimensionalized version  (\ref{forward_adim2}) of (\ref{balance_escalar}), see 
Appendix A for details and typical parameter choices.


\subsection{Approximate observation operators and discrete Bayesian formulation}
\label{sec:observation_aprox}

For computational studies, the original problem is replaced by a discretized version 
that employs approximate observation operators $\mathbf o_{\rm ap}$
\begin{eqnarray} \begin{array}{ll}
\mathbf o_{\rm ap}: & \mathbb R^P \longrightarrow \mathbb R^D 
\\ [1ex]  
&\nnu \longrightarrow (u_{\nnu, {\rm ap}}(\mathbf r_j, t_m))_{j=1,\ldots,J, m=1,\ldots,M},
\end{array} \label{observation_ap} 
\end{eqnarray}
where $u_{\nnu, {\rm ap}}$ is the approximation to the solution of (\ref{balance_half})
generated by solving numerically  (\ref{forward_adim2}) with the selected schemes.

Here, we will use the scheme (\ref{discretization}) described in Appendix C with three
different mesh choices: adapted, stratified and uniform. 
 The latter two lead to two different and fixed approximate 
observation operators during the whole optimization and sampling procedures, which 
we denote $\mathbf o_{ap_s}$, and $\mathbf o_{ap_u}$, respectively.
The posterior distribution to be explored becomes
\begin{eqnarray}
p(\nnu |  \dn ) \sim \exp \left( -\frac{1}{2} \| \mathbf o_{\rm ap}( \nnu) - \dn  \|^2_{\Gn^{-1}}
-\frac{1}{2} \|\nnu - \nnu_0\|_{\Gpr^{-1}}^2 \right),
\label{unnormalized_ap}
\end{eqnarray} 
whereas the cost to be optimized to evaluate the MAP estimate is
\begin{eqnarray}
J(\nnu) = \frac{1}{2} \| \mathbf o_{\rm ap}( \nnu) - \dn  \|^2_{\Gn^{-1}}
+ \frac{1}{2} \|\nnu - \nnu_0\|_{\Gpr^{-1}}^2.
\label{bcost_ap}
\end{eqnarray}
Compared to the original distributions (\ref{unnormalized}) and costs (\ref{bcost}) 
the difference $\mathbf o(\nnu) - \mathbf o_{\rm ap}(\nnu)$ could result in variations
in the original landscape of high probability configurations (resp. local minima),
specially when the mesh is kept fixed.
Instead, the use of meshes that adapt to inclusions as they vary results in continuous 
changes in the approximate observation operator in an attempt to reproduce better 
the original observation operator.


Next, we introduce techniques to numerically approximate the MAP points.


\section{Calculation of maximum a posteriori estimates}
\label{sec:optimization}

To approximate the highest probability parameter sets which minimize
the cost functionals (\ref{bcost}) or (\ref{bcost_ap}), we adapt 
techniques  of deterministic optimization \cite{nocedal} (more precisely, 
Levenberg-Marquardt-Fletcher type approaches \cite{lmf}) and also
Markov Chain Monte Carlo methods.
A delicate point in the implementation of optimization schemes is 
the approximation of the derivatives of solutions $u_{\nnu}$ of 
(\ref{balance_escalar}) with respect to $\nnu$, 
which requires either a characterization of such derivatives as the
solution of an auxiliary problem or adequate numerical schemes
depending on the selected type of mesh. 
We introduce next an automatic strategy to minimize functional
(\ref{bcost}) employing adaptive meshes and an strategy to locate
possible additional minima resorting to fixed meshes.
In principle, we can locate additional minima by means of 
affine invariant MCMC samplers \cite{goodman} too.
MCMC methods are technically easier to  implement, since only forward 
solvers are needed. However, the number of auxiliary problems to be 
solved increases from a few tens to a few hundred thousand or millions.


\subsection{Iterative optimization scheme for adaptive meshes}
\label{sec:oadaptive}

We choose the prior mean $\nnu_{0}$ as initial guess of the parametrization, 
$\nnu^0 = \nnu_{0}$ and implement the Newton type iteration  
$ \nnu^{k+1} = \nnu^{k} + \boldsymbol \xi^{k+1}$  where  
$ \boldsymbol \xi^{k+1}$ is the solution of  
\begin{eqnarray}
\left(\mathbf H(\nnu^k) + \omega_k
{\rm diag}(\mathbf H(\nnu^k) ) \right)
\boldsymbol \xi^{k+1} = -  \mathbf g(\nnu^k)
\label{it_false}
\end{eqnarray}
Here, $\mathbf H(\nnu)$ and $\mathbf g(\nnu)$ represent the Hessian and 
the gradient of the selected cost, that involve first and second order
derivatives of the solutions of the forward problems with respect to $\nnu$.
Since we are only interested in seeking descent directions, in practice, we 
replace  $\mathbf H(\nnu)$ by the Gauss-Newton  part of the  Hessian 
$\mathbf H^{\rm GN}(\nnu)$ neglecting  second order derivatives.
The small parameter $\omega_k >0$ is adjusted to guarantee a decrease in 
the cost.  We reduce it as the cost $J(\nnu^k)$ decreases and we increase it 
if the proposed  $\boldsymbol \xi^{k+1}$ leads to a gain in the cost \cite{lm,lmf}.

In the adaptive approach we propose, we recalculate the mesh and triangulation 
of the computational region $R$ any time we update the inclusion parametrization
$\nnu^k$ in (\ref{discretization})-(\ref{matrices_artificial}) with $\rho$, $\chi$, 
$v_{\rm p}$ defined in (\ref{coefrho})-(\ref{coefs}) and $\gamma = \chi/v_{\rm p}$.
We adapt the finite element mesh to the shape  of the domains in which the 
physical properties take different constant values using a technique developed 
by Gilbert Strang and Per-Olof Persson in  \cite{Persson,Strang-Persson}. Then, 
$u_{\nnu^k}$ is calculated by means of (\ref{discretization})-(\ref{matrices_artificial}). 
At each step, we approximate the derivatives of the solutions $u_{\nnu^k}$ with 
respect to $\nnu$  by means of algorithmic differentiation  \cite{ForwardDiff}.
The process stops when either the difference between parameterizations or the
cost value fall below given tolerances. The final value provides an approximation
to $\nnumap$.

The algorithm steps can be summarized as follows, see \cite{oleaga} for implementation
details:
\begin{itemize}
\item  {\it Initialization:} Define 
\begin{itemize}
\item prior mean $\nnu_0$ (inclusion shape and material parameters) and prior covariance  $\Gpr$,
\item measured data $\dn$ and noise level $\sigma_{\rm noise}$,
\item emitter grid $\mathbf x_k$, $k = 1, \ldots, K$, receiver grid $\mathbf r_j$, 
$j = 1, \ldots, J$, recording times $t_m,$ $m=0,\ldots,M$,
\item material parameters for $\rho$ and $v_{\rm p}$ on the layered structure and emitted signal $f(t)G(\mathbf x)$,
\item maximum number of optimization steps $S$, spatio-temporal steps $\delta t$, $\delta x$  and $\delta y$ for the solver, final time $T$ and  tolerances.
\end{itemize}
Set the initial parameterization $\nnu^0$ equal to the prior mean. 

\item {\it Optimization:}  From $i=1$ to $S$
\begin{itemize}
\item Build a triangulation ${\cal T}^i$ adapted to the current geometry of the inclusions and the underlying  layered structure (ie. each triangle is fully contained in a subdomain with constant coefficients):
\begin{itemize}
\item Implement the Persson-Strang method \cite{Persson,Strang-Persson}.
\item Store the resulting points and triangulation.
\end{itemize}

\item Solve numerically the boundary value problem (\ref{forward_adim2}):
\begin{itemize}
\item Construct a space of linear finite elements on the adapted triangulation ${\cal T}^i$, designed to admit forward automatic differentiation with respect to points, triangulation and elastic parameters.
\item Update the FEM matrices (\ref{matrices_artificial}) to reflect the coefficient values for the current inclusion.
\item Implement the variant of the discretization scheme (\ref{discretization}) for (\ref{2nd}).
\item Evaluate the numerical solution $u_{\nnu^i}$  at the receivers to calculate the observation operator $\mathbf o_{\rm ap}(\nnu^i)$.
\end{itemize}

\item Solve system (\ref{it_false}) to calculate $\nnu^{i+1}$:
\begin{itemize}
\item Approximate the variation of the solutions $u_{\nnu^i}$ at the receivers with 
respect to $\nnu$ by means of algorithmic differentiation  \cite{ForwardDiff} to obtain
$\mathbf H(\nnu^i)$ and $\mathbf g(\nnu^i)$.
\item Propose a descent direction $\boldsymbol \xi^{i+1}$ using (\ref{it_false}) and 
evaluate the cost $J(\nnu^{i+1})$ given by (\ref{bcost_ap}) with 
$\nnu^{i+1} = \nnu^{i} + \boldsymbol \xi^{i+1}$.
\item If $J(\nnu^{i+1})< J(\nnu^{i})$ accept $\nnu^{i+1}$, set $\omega^{i+1}=
\omega^{i+1}/2$ and move to the next optimization step. Otherwise, divide 
$\omega^i$ by 2 and repeat until satisfied.
\item If either $\| \nnu^{i+1} - \nnu^{i}\|$ or $J(\nnu^{i+1})$ fall below specified tolerances,
stop.
\end{itemize}
\end{itemize}

\item {\it Output:} Optimal inclusion shape and material parameters defining $\nnumap$, 
gradient and Hessian of the cost at  $\nnumap$, as well as intermediate $\nnu^i$ and
$J(\nnu^i)$ values.
\end{itemize}

This algorithm can be coded in Julia to exploit specific meshing and automatic 
differentiation packages. The resulting adaptive scheme is fully automatic: the 
user must only provide the data and the parameters defining the cost functional, 
that is, the prior knowledge, and the level of noise.
Meshes are automatically adapted to the new object and internal variables
for the calculation of derivatives and the advance of iterations are adjusted
automatically. The user must only provide the data and the parameters defining
the cost functional, that is, the prior knowledge, and the level of noise.  

On a laptop, building an adapted triangulation takes about $31$ seconds and solving 
one forward problem takes about $29$ seconds. One iteration of the optimization routine 
takes about $14$ minutes (it may involve several object proposals and triangulations). 
The usual $20$ steps make take about 4 hours. This process
is not paralellized due to the difficulty of parallelizing mesh generation.

\begin{figure}[!hbt]
\centering
\includegraphics[width=12cm]{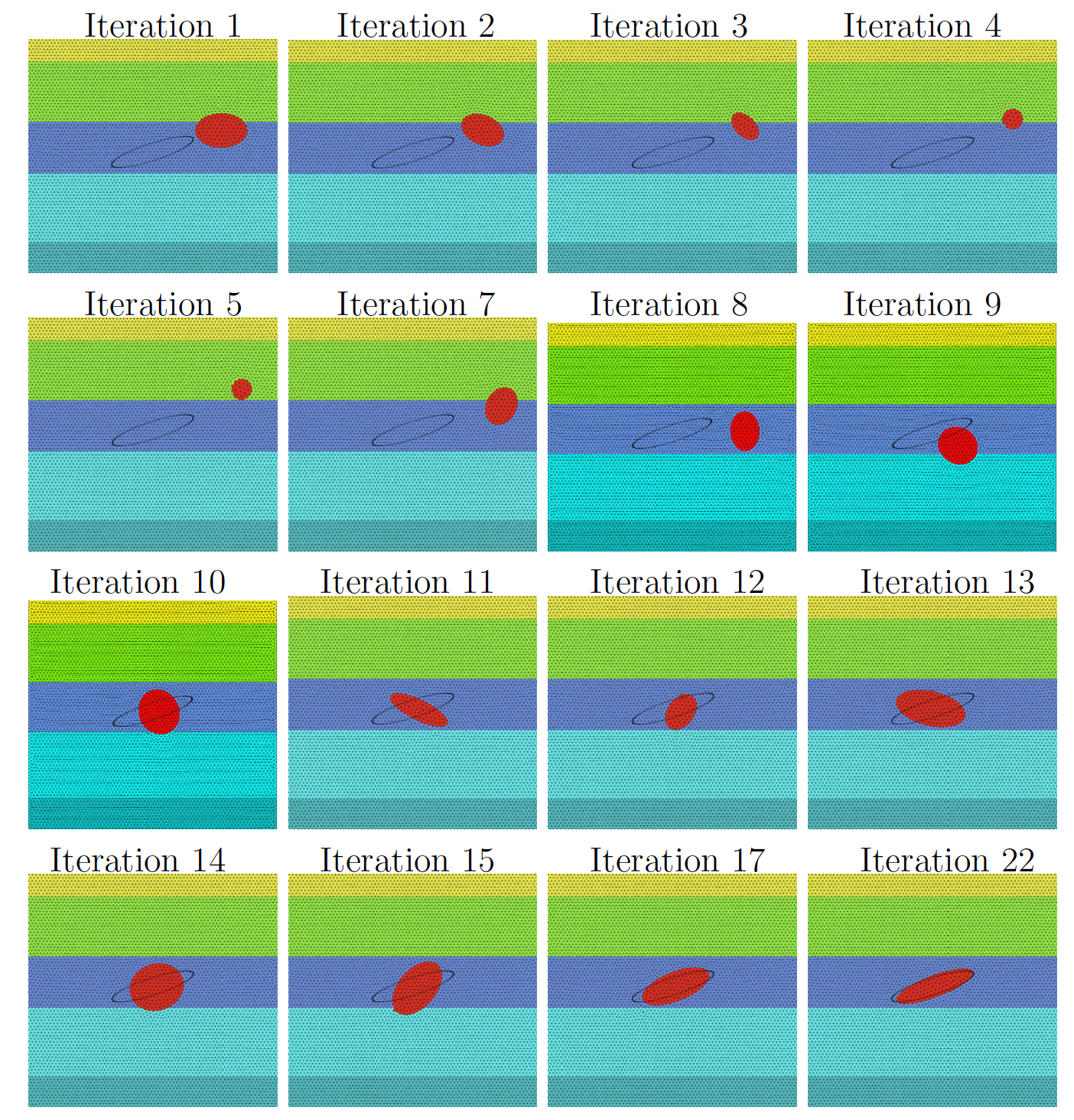}   
\caption{True inclusion (black curve)  compared to  the shapes obtained at successive
iterations (red) for $10\%$ noise during adaptive constrained optimization.}
\label{fig4}
\end{figure}

\begin{figure}[!hbt]
\centering
\includegraphics[width=12cm]{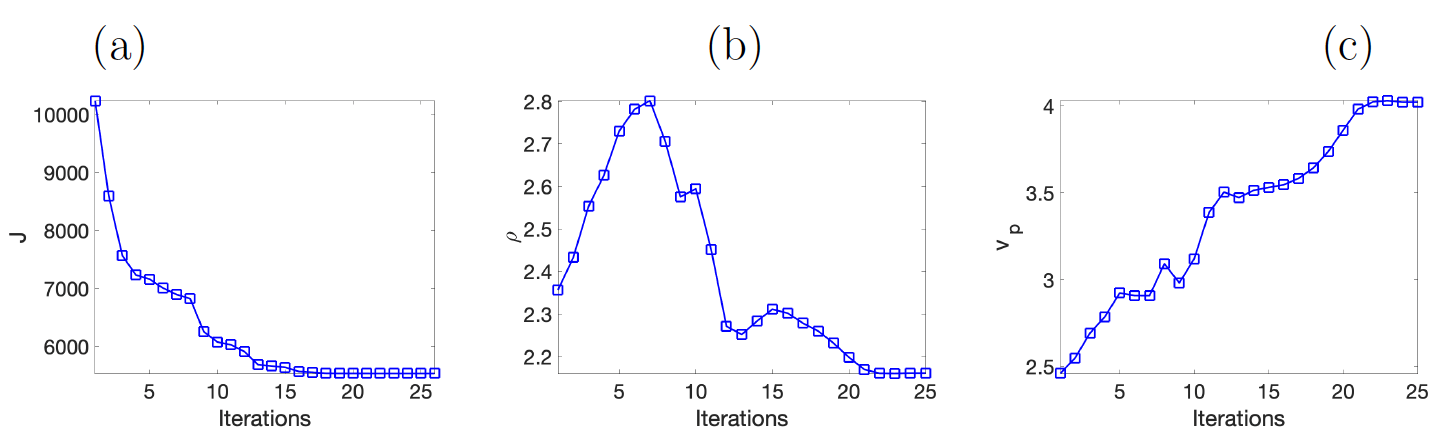} 
\caption{Evolution of  (a) the cost,  (b) $\rho$ and (c) $v_{\rm p}$ inside the inclusion
along the iterations for the simulation in Figure \ref{fig5}.}
\label{fig5}
\end{figure}

\begin{figure}[!hbt]
\centering
\includegraphics[width=12cm]{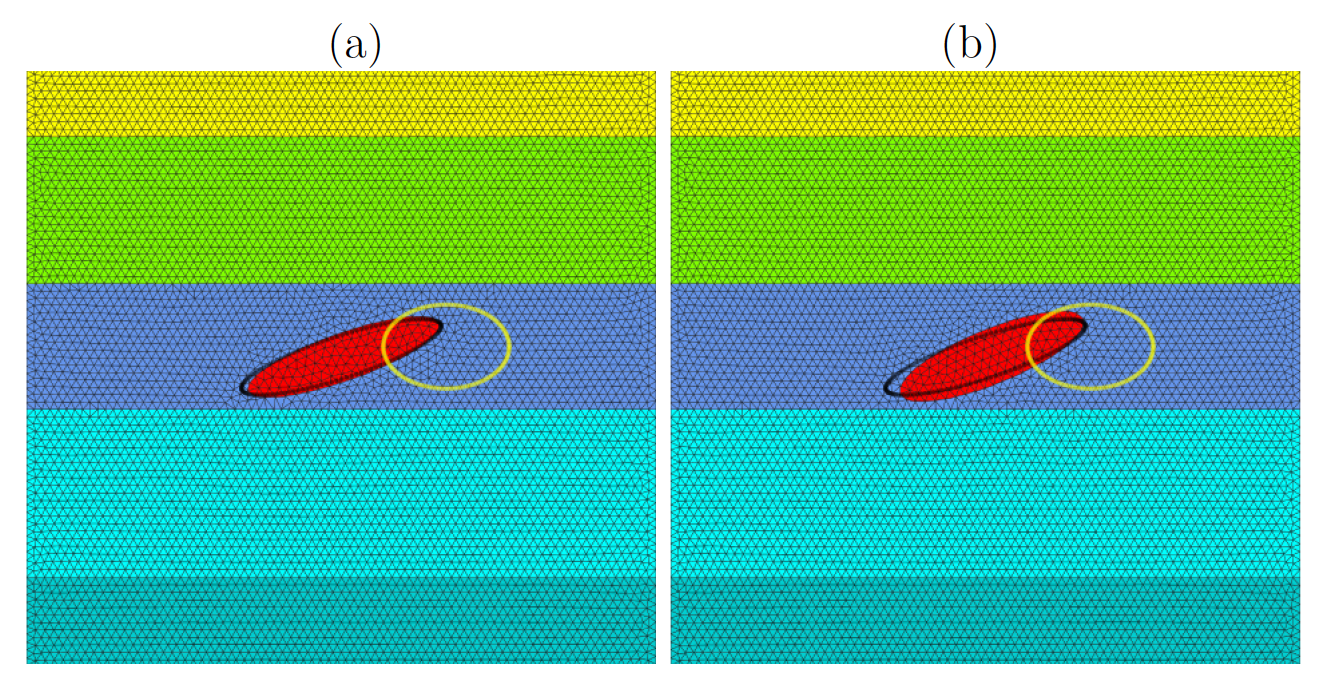}  
\caption{True inclusion (black curve) compared to the initial (prior) configuration (yellow
curve) and the MAP estimate of the shape (meshed in red) after convergence of 
adaptive constrained optimization in (a) $28$  iterations for $5\%$ noise, (b)  $29$ 
iterations for $15\%$ noise. MAP values are given in Table \ref{table2}.}
\label{fig6}
\end{figure}

Figures \ref{fig4}-\ref{fig6} summarize results for the true configuration 
represented in Figure \ref{fig1} for increasing noise levels in the data. 
Synthetic data for different noise levels $r$ are generated an explained in 
Appendix A. 
Figures \ref{fig4} and \ref{fig5} illustrate the evolution of the inclusion geometry
and its material properties along the optimization procedure for $10\%$ noise. 
The object can move from one layer to another, shrink and expand or rotate 
as its material properties change.
Table \ref{table2} contains the parameter values for the true inclusion and 
the  initial approximation, as well as the MAP estimate obtained optimizing 
(\ref{bcost}) for $5\%$ and $15\%$ noise. Figure \ref{fig4} visualizes the
true object, the prior and the final estimate provided by the optimization
procedure. Inclusions are defined by the  parameters  
$\nnu = (c_{\rm x},c_{\rm y}, a,b,\theta,\rho,v_{\rm p})$, that its, 
center, semi-axes, angle with the positive $x$-axis, density and wave speed.
We choose as prior means for $\rho$ and $v_{\rm p}$  the average of the 
maximum and minimum values available for the layers, that is, 
$(\rho_{\rm max} + \rho_{\rm min})/2$ 
and $(v_{{\rm p}, {\rm max}} + v_{{\rm p}, {\rm min}})/2 $.
Then, the  standard deviations $\sigma_\rho$ and $\sigma_{v_p}$ are taken 
equal to half the difference between the maximum and minimum values in 
the whole sequence of known layers, that is,
$ \sigma_\rho = (\rho_{\rm max} - \rho_{\rm min})/2$ and
$\sigma_{v_{\rm p}} = (v_{{\rm p}, {\rm max}} - v_{{\rm p}, {\rm min}})/2$.
For the tests  presented here we have set  in the cost  
$\Gpr =$  diag $(1$, $1$, $0.5$,  $0.5$,  $0.1$, $\sigma_\rho^2$, 
$\sigma_{v_{\rm p}}^2$) and $\Gn= \sigma_{\rm noise}^2 \mathbf I$,
$\sigma_{\rm noise} = \sigma \, r/100$,   $\sigma$ being the normalized
$\ell^2$ norm of the true data (see Appendix A).

\begin{figure}[!hbt]
\centering
\includegraphics[width=12cm]{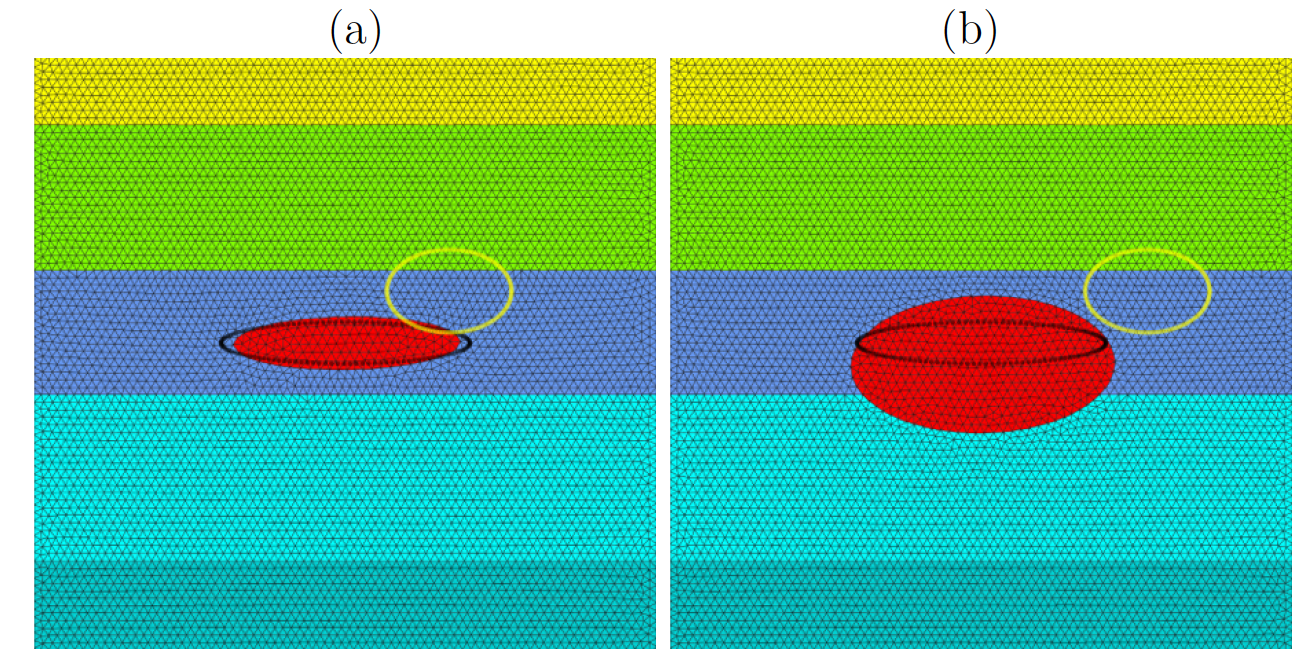}  
\caption{MAP estimates of the shape of an horizontal inclusion obtained by adaptive 
constrained optimization with slightly different prior shapes for $10\%$ noise.
The MAP estimates of $\rho$ and $v_{\rm p}$ are $2.13$ and $4.03$ for (a) and 
$2.32$ and $3.96$ for (b).}
\label{fig7}
\end{figure}

In the tests we have performed, this automatic algorithm leads to similar 
results for rotated inclusions (clockwise or anti-clockwise) even if we place the 
prior in a different layer or if the true object lies between different layers.  
Keeping the same material parameters, the MAP estimates for the geometry
of  rotated inclusions in this set-up remain quite close to the true shape, as well 
as the density, while the velocity varies a bit more.
Remarkably, the prediction for horizontal inclusions may depart from the true 
shape as we vary the prior, see Figure \ref{fig7}, obtained for different initial curves
with the same material parameters.
Notice that changing the prior we change the continuous functional to be 
optimized. Later MCMC studies will clarify this observation, see Figure \ref{fig8}(d).

\vspace{0.2cm}
\begin{table}[hbt!]
\begin{tabular}{cccccc}
  \hline
  \textbf{Parameter} & \textbf{True} & \textbf{Prior} & \textbf{Adaptive $5$\% } 
  & \textbf{Adaptive $15$\% }\\\hline
  x center & 0.0 & 0.5 & 0.0110207 & 0.030625 \\
  y center & -1.45 & -1.4 & -1.44922 & -1.446588 \\
  semi-major axis & 0.5 & 0.3 & 0.485523 & 0.472986 \\
  semi-minor axis & 0.1 & 0.2 & 0.113252 & 0.1305437 \\
  rotation angle & 0.314159 & 0.0 & 0.340005 & 0.38868 \\
  density & 2.1 & 2.3 & 2.13201 & 2.182858 \\
  velocity & 4.4 & 2.4 & 4.17327 & 3.9267512\\ \hline
  \end{tabular}
\caption{Dimensionless inclusion parameters compared to the parameters defining 
the initial approximation and the final MAP estimate obtained by adaptive constrained 
optimization for different noise levels. 
} 
\label{table2}
\end{table}

\subsection{Iterative optimization scheme for fixed meshes}
\label{sec:ofixed}

When working with prefixed meshes, we implement the previous scheme with 
four differences.

First, we construct the mesh at the start, and keep it unchanged during the 
optimization process. The observation operator is thus (\ref{observation_ap})
with $\mathbf o_{\rm ap_s}$ or $\mathbf o_{\rm ap_u}$, and the cost 
(\ref{bcost_ap}).
Thanks to the structure of the scheme (\ref{discretization}), all the matrices can be
calculated for a configuration without inclusions and stored at the start of the 
process. Whenever we propose an inclusion, only the coefficients corresponding 
to nodes involving it are updated. The resulting matrices are kept fixed
to calculate the approximate solution $u_{\nnu^k, {\rm ap}}$ of the corresponding 
problem and the required observed data. Matlab routines to construct FEM 
matrices allow to locate the affected nodes, enabling us to perform this process 
quite fast. 

Second, to estimate the change of the solution $u_{\nnu ,{\rm ap}}$ when varying 
the parameters $\nu_i$, we rely on difference approximations \cite{elena}. 
\[
{\partial u_{\nnu, {\rm ap}} \over \partial \nu_i}(\mathbf r_k,t_m) 
\sim {u_{\nnu+\eta_i, {\rm ap}}(\mathbf r_k,t_m) - 
u_{\nnu, {\rm ap}}(r_k,0,t_m)  \over \eta_i},
\]
with $\eta_i$ small, 
$u_{\nnu+\eta_i, {\rm ap}}$ being the  numerical solution of the forward
problem with $\nu_i$ replaced by $\nu_i+\eta_i.$ Notice that we do not need
detailed derivative studies. We just use these quotients to generate descent 
directions for the cost.

Third, to calibrate the steps $\eta_i$ in these approximations we have 
swept a grid of  small values, that are selected and then kept fixed for all 
approximations during the iterations performed in the optimization 
procedure.

Four, this strategy results in iterative procedures that decrease the initial value
of the cost and stop when changes in the cost or parameters fall below
given tolerances. While varying the parameter $\eta$, the scheme may converge 
to different parameter sets $\nnu$ with different cost values.  The set yielding the 
minimum value of the cost is selected as $\nnumap$.

For the data considered in Section \ref{sec:oadaptive}, Table \ref{table3} 
collects the results obtained for $\nnumap$ employing fixed uniform meshes
with observation operator $\mathbf o_{{\rm ap}_u}$.
Similar results are found for $\mathbf o_{{\rm ap}_s}$. In both cases,
the optimization algorithm provides other configurations with slightly higher
cost as we sweep the grid for $\eta$, which suggests the presence of
secondary local minima, since the prior and the mesh are kept fixed now.
Nevertheless, most choices of $\eta$ lead to $\nnumap$.
Compared to the adaptive algorithm, carrying out this process on a laptop for 
a single value of $\eta$  takes a few minutes, when coded in Matlab and 
conveniently parallelized. For the grids we are considering, with maximum
steps $0.04$ in $x$ and $y$, it takes about $1.5$-$3$ minutes to converge 
for a uniform mesh and  between $14$ and $94$  minutes for a stratified mesh. 
Switching to maximum step $0.06$ in the $x$ direction, time goes down to 
$4$-$14$ minutes for the stratified mesh.
Sweeping a grid of values of $\eta$ to search for  other possible local minima 
increases the cost: about $20$ minutes for the uniform mesh and about $5$
hours or $1$ hour for the stratified mesh. This requires checking local convergence 
to asses whether the stopping point is a possible local minimum or the iterations 
stops because the descent strategy does not produce suitable descent directions. 
As said before, the adaptive algorithm  is more difficult to parallelize and it may
take a few hours on a laptop without parallelization. However, it is fully automatic. 

\begin{table}[hbt!]
\begin{tabular}{cccccc}
  \hline
  \textbf{Parameter} & \textbf{Fixed $5$\%} & \textbf{Fixed $15$\%} & 
  \textbf{MCMC $5$\% } & \textbf{MCMC $15$\% }\\ \hline
  x center             & 0.01766   &  0.040319 & 0.0217    & 0.0052    \\
  y center             & -1.4489   &  -1.448       & -1.4476   & -1.4548    \\
  semi-major axis & 0.4844    & 0.4766       & 0.4856     & 0.4588    \\
  semi-minor axis & 0.10071  &  0.13516     & 0.0940     & 0.1379   \\
  rotation angle    & 0.3264    &  0.3474       & 0.3334     & 0.3175    \\
  density               & 2.0977    & 2.1991       & 2.0718     & 2.1457   \\
  velocity              & 4.3678     & 3.8257       & 4.4645     & 3.9275    \\ \hline
  \end{tabular}
\caption{Counterpart of Table 3 for constrained optimization and MCMC 
schemes using fixed uniform meshes.} 
\label{table3}
\end{table}


\subsection{Optimization by Markov Chain Monte Carlo methods}
\label{sec:omcmc}

A  Markov chain is a sequence of events for which the probability of 
an event depends only on the event just before it. To define
a Markov Chain we need three elements \cite{mcmc}:
the space of states  $X$, that is, the set of values the chain can take,
the transition operator  $q(x^{k+1}|x^k)$, which defines the probability
of transitioning from state $x^k$ to state $x^{k+1}$, and
the initial distribution  $\pi$, which defines the probability of being in 
any one of the possible states at the start, for $k=0$.
To generate a Markov chain $x^0 \rightarrow x^1 \rightarrow x^2 
\ldots \rightarrow x^k\rightarrow  \ldots$ we sample the initial state 
$x^0$ from $\pi$ and transition from $x^k$ to $x^{k+1}$ according to
$q(x^{k+1}|x^k)$, $k \geq0$. 

Markov chains satisfying a number of properties (time homogeneity,
detailed balance) are shown to equilibrate to target distributions $p$
under some conditions \cite{inference,mcmc}.
Different strategies to construct Markov Chains enjoying these properties
have been proposed. We choose here the affine invariant ensemble sampler   
developed in \cite{goodman}, that can handle multimodal distributions and 
allows for parallelization \cite{sergei}.
The idea is to create $W$ chains that are mixed at each step:
\cite{goodman}:
\begin{itemize}
\item Set the number of chains $W$, the number of steps $S$
         and choose $a  \!\sim\! 2$.
\item Set the sample space $X=\mathbb R^P$, $P$ dimension
          of the parameter space.
\item Draw $\nnu^0_1,\ldots,  \nnu^0_W$ in $X$
         with probability $p_{\rm pr}$.     
\item From $k=0$ to $k=S$:
         \begin{itemize}
            \item Generate a permutation $\sigma$ of $(1,\ldots,P)$ without 
              fixed elements.
            \item For $w= 1, \ldots, W$
            \begin{itemize}
                \item  Draw $ z_w $ from $g(s) = s^{-1/2}$ 
                if $s \in [a^{-1},a]$, $0$ otherwise.
                \item Set $\nnu_{{\rm prop},w}= \nnu_{\sigma(w)}^k
               + z_w (\nnu_{w}^k-\nnu_{\sigma(w)}^k)$.
                \item  Calculate the acceptance rate
               $\alpha={\rm min}\left(1, z_w^{P-1}{ p_{\rm pt}(\nnu_{{\rm prop},w})
                \over p_{\rm pt}(\nnu^k_w) } \right)$.
                \item  Draw $u$ with uniform probability $U(0,1)$.
                If $u < \alpha$, $\nnu^{k+1}_w = \nnu_{{\rm prop},w}$,
                otherwise $\nnu^{k+1}_w = \nnu^k_w$.
           \end{itemize}    
         \end{itemize}
\end{itemize} 
Discarding an initial transient stage formed by $B$ samples on account of chain 
equilibration, the remaining states $\nnu^k_w$ the chain takes  sample from the 
probability $p_{\rm pt}$.  Notice that we can sample unnormalized distribution like
(\ref{unnormalized}) because the normalization factors scale out.
This algorithm needs $W>2P$ to properly sample the target posterior distribution. 

Due to the high computational cost of recalculating meshes and auxiliary matrices
for each parameter proposal, we employ this method  to sample $p_{\rm pt, ap}$ 
given by (\ref{unnormalized_ap}) for observation operators $\mathbf o_{\rm ap}$ defined
on fixed meshes, that is, $\mathbf o_{\rm ap_s}$ and $\mathbf o_{\rm ap_u}$. 
The MAP estimate $\nnumap$ is the sample for which the probability is higher, that is, 
the cost $J(\nnu) = -\log(p_{\rm pt, ap}((\nnu) ))$ is minimum.

\begin{figure}[!hbt] \centering
\includegraphics[width=12cm]{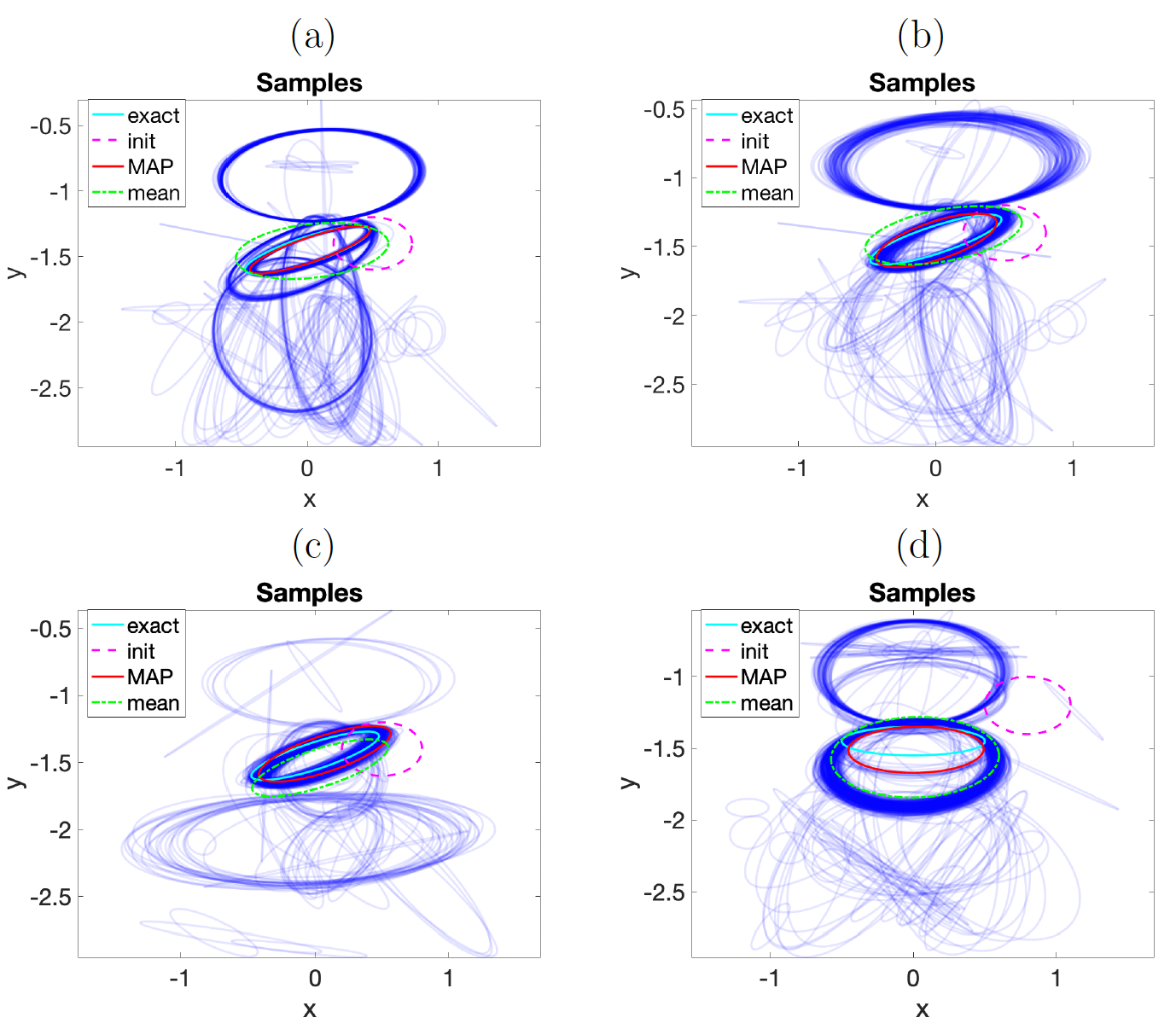}
\caption{Last $W$ samples for (a) $5\%$ noise, (b) $15\%$ noise, (c) $15\%$ noise
with a bottom layer of density $\rho=2.1$ and larger velocity $v_{\rm p}=4.4$, (d) $10\%$
noise with horizontal inclusion. Parameter values $W=480$, $a=2$, $K=1000$.}
\label{fig8}
\end{figure}

Table \ref{table3} collects the parameters obtained for $\nnumap$ for the test case
in Figure \ref{fig1} employing  fixed uniform meshes and keeping the data used in
the previous two sections. The MAP estimates are similar in the three cases.
Figure \ref{fig8} illustrates additional high probability configurations obtained in this 
way, as suggested by the sample concentration. Panel (c) suggests that the increased
signal to noise ratio associated to a more reflective bottom layer reduces uncertainty.
We discuss this point further in the next section.

\section{Uncertainty quantification in the object based full-waveform inversion problem}
\label{sec:uncertainty}

Maximum a posteriori estimates provide a prediction of the most likely values for the 
inclusion. More precise uncertainty studies provide estimates of the expected range
of variation of the different parameters and derived magnitudes.


\subsection{Uncertainty quantification based on the Laplace approximation}
\label{sec:ulaplace}

The Laplace approximation linearizes (\ref{unnormalized_ap}) about $\nnumap$ and
approximates it by a multivariate Gaussian distribution  ${\cal N}(\nnumap, \Gpo)$
The approximate posterior covariance is  
\begin{eqnarray*}
\Gpo =  \mathbf H(\nnumap)^{-1} \sim ( \mathbf F(\nnumap)^t \Gn^{-1} 
\mathbf F(\nnumap) + \Gpr^{-1})^{-1} = \mathbf H^{\rm GN}(\nnumap)^{-1},
\end{eqnarray*}
where  $\mathbf H^{\rm GN}(\nnumap)$ represents the Gauss-Newton
part of the Hessian $ \mathbf H(\nnumap)$, which neglects second order derivatives, 
and  $\mathbf F(\nnumap)= \left(\frac{\partial u_{\nnumap}}{\partial \nu_i}(p_k) 
\right)_{k,i}$  with
\[ \mathbf p = ((\mathbf r_1,t_1), \ldots, (\mathbf r_J,t_1), \ldots,
(\mathbf r_1,t_M), \ldots, (\mathbf r_J,t_M)). \]
Samples of the Gaussian approximation ${\cal N}(\nnumap, \Gpo)$ are generated 
as
\begin{eqnarray}
\nnu = \nnumap+ \Gpo^{1/2} 
\mathbf w, \label{sample_map}
\end{eqnarray}
$\mathbf w$ being a standard iid vector.  

\begin{figure}[hbt!]
\centering
\includegraphics[width=12cm]{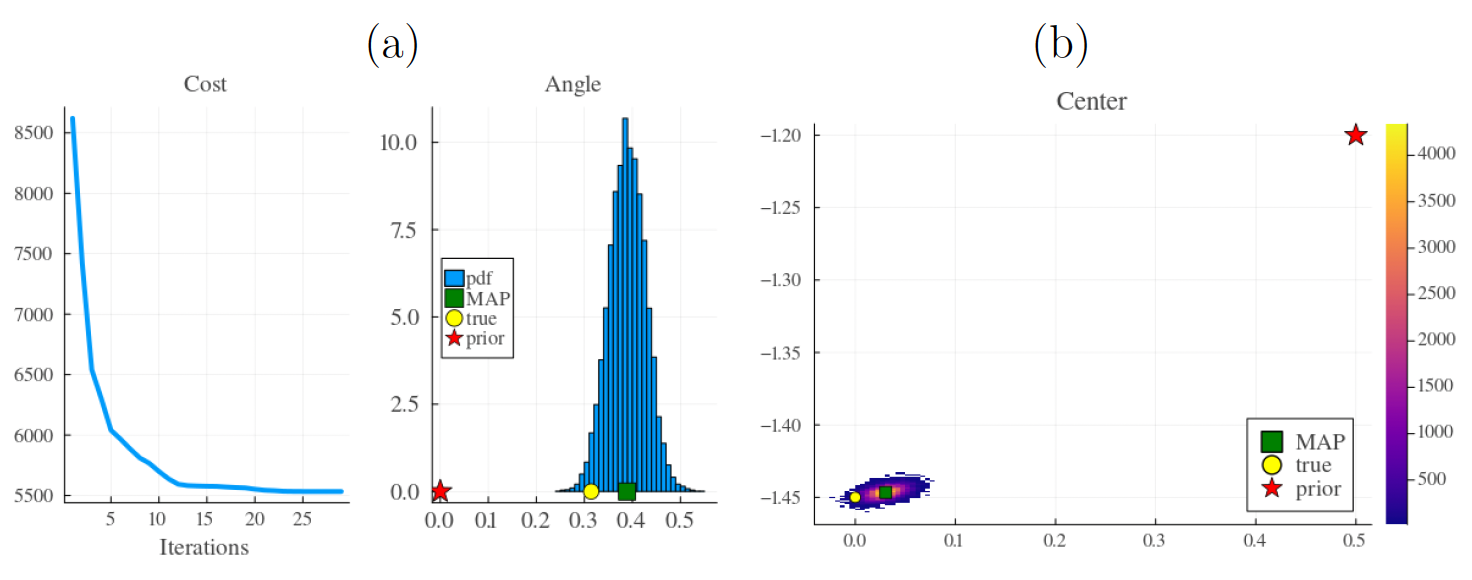} 
\caption{For Figure \ref{fig6}(b) with 15\% noise:
(a) Cost evolution during MAP calculation by adaptive optimization
and uncertainty in angle orientation.
(b) Uncertainty in the center position. }
\label{fig9}
\end{figure}

\begin{figure}[!hbt]
\centering
\includegraphics[width=12cm]{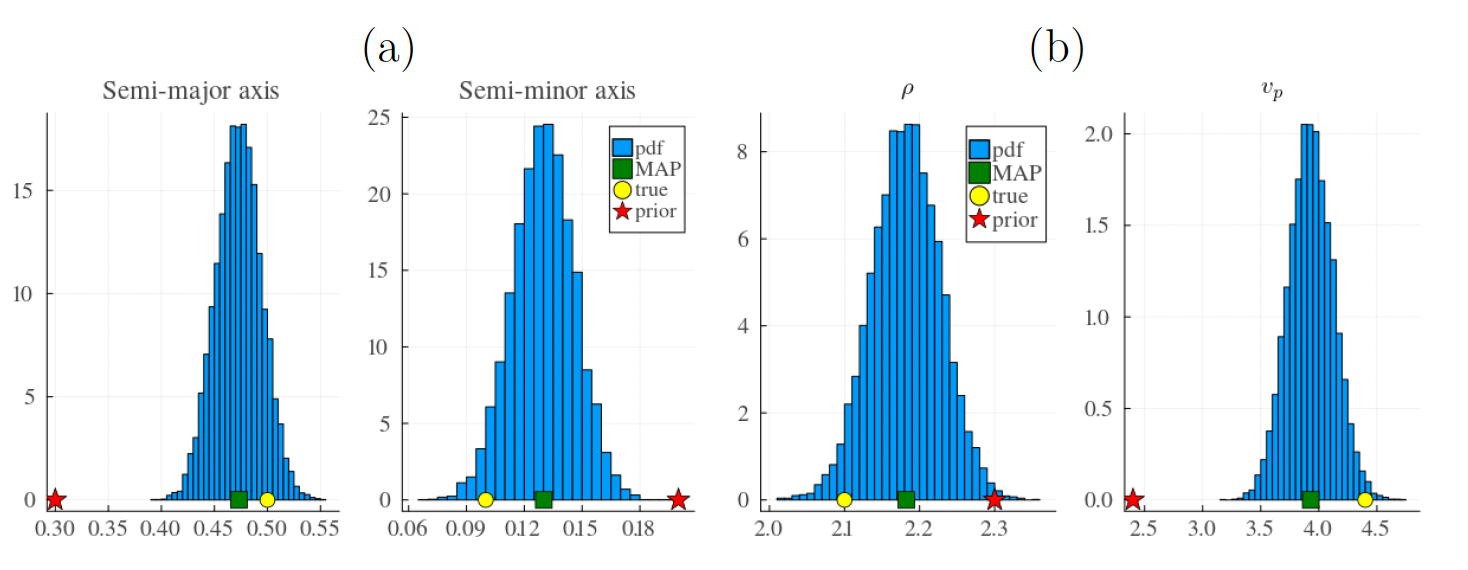} 
\caption{For Figure \ref{fig6}(b) with 15\% noise:
Uncertainty in (a) the semi axes estimate and (b)
the material parameters.}
\label{fig10}
\end{figure}

Figures \ref{fig9}-\ref{fig10} visualize the information we obtain about the 
uncertainty ranges in our predictions for the test case considered in Section
\ref{sec:oadaptive}, based on $10.000$ samples.
We compare the true values, the initial values, the MAP estimate and
superimpose histograms built sampling the Laplace approximation of
a the posterior distribution (\ref{sample_map}) with a pdf normalization.
These figures use the outcome of the automatic calculations performed
with the adaptive optimization scheme in Section \ref{sec:oadaptive}.

The scheme in Section \ref{sec:ofixed} produces similar results for 
$\nnumap$, but suggests the presence of additional local minima
representing additional high probability configurations.
The posterior probability for the discretized observation operator
(\ref{observation_ap}) might be a multimodal distribution with 
additional peaks. Notice that this approach uses an approximation of 
the observation operator based on fixed meshes, see remarks in 
Appendix C, keeping the same data.   


\subsection{Uncertainty quantification based on MCMC studies}
\label{sec:umcmc}

More detailed information on uncertainty on parameter ranges and derived
magnitudes  is extracted from the analysis of the samples provided by the 
MCMC chains.

\begin{figure}[!hbt]
\centering
\includegraphics[width=12cm]{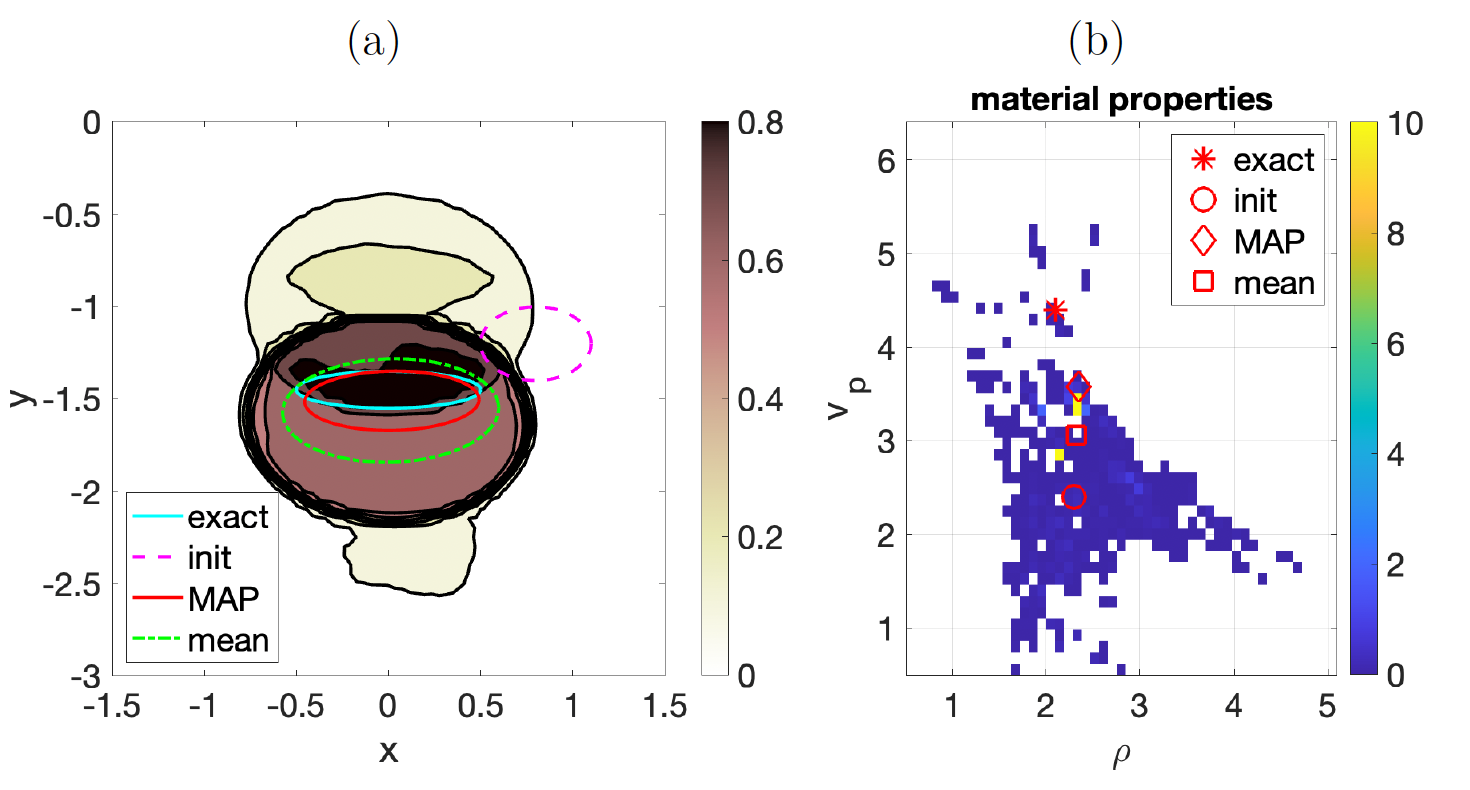} 
\caption{(a) Contour plot of the probability of belonging to the inclusion 
and (b) projection of the histogram for its material properties 
(pdf normalization) calculated by MCMC sampling and 
corresponding to Figure \ref{fig8}(d) with $10\%$ noise.}
\label{fig11}
\end{figure}

\begin{figure}[!hbt]
\centering
\centering
\includegraphics[width=12cm]{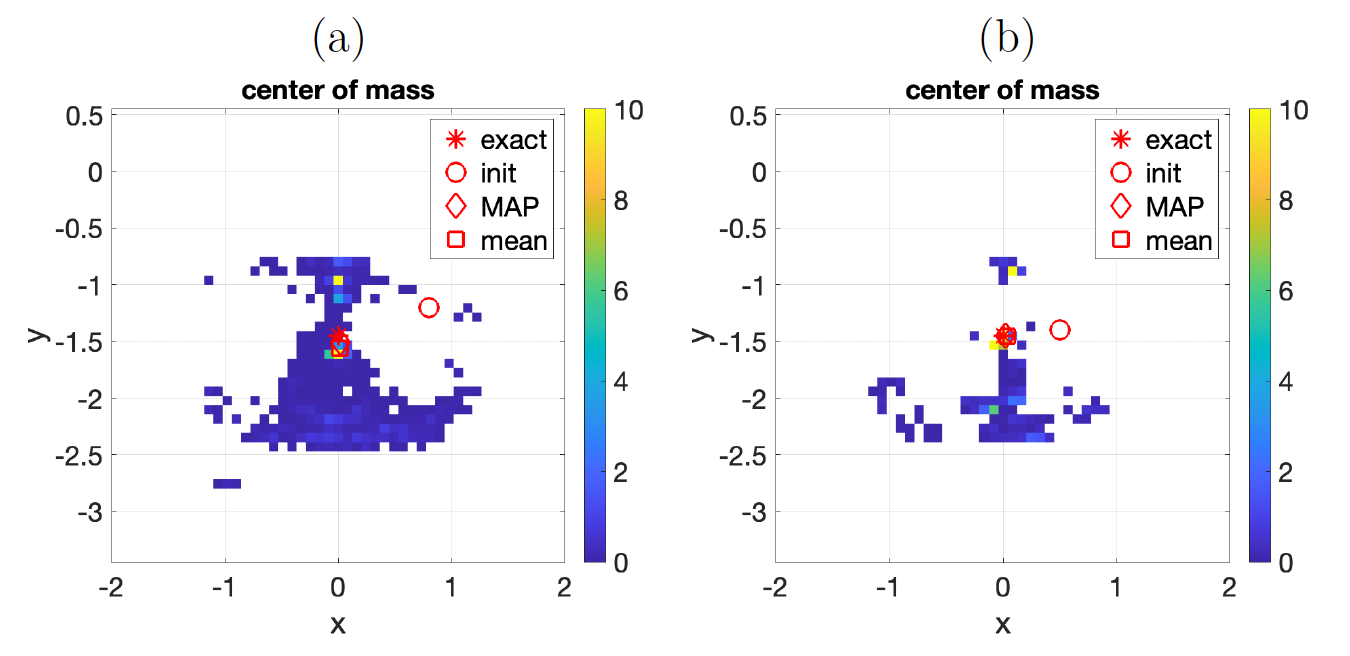} 
\caption{Histograms (pdf normalization) for the centers of mass calculated 
by MCMC sampling: (a) corresponding to Figure \ref{fig11},
(b) corresponding to Figure \ref{fig13}.}
\label{fig12}
\end{figure}

\begin{figure}[!hbt]
\centering
\includegraphics[width=12cm]{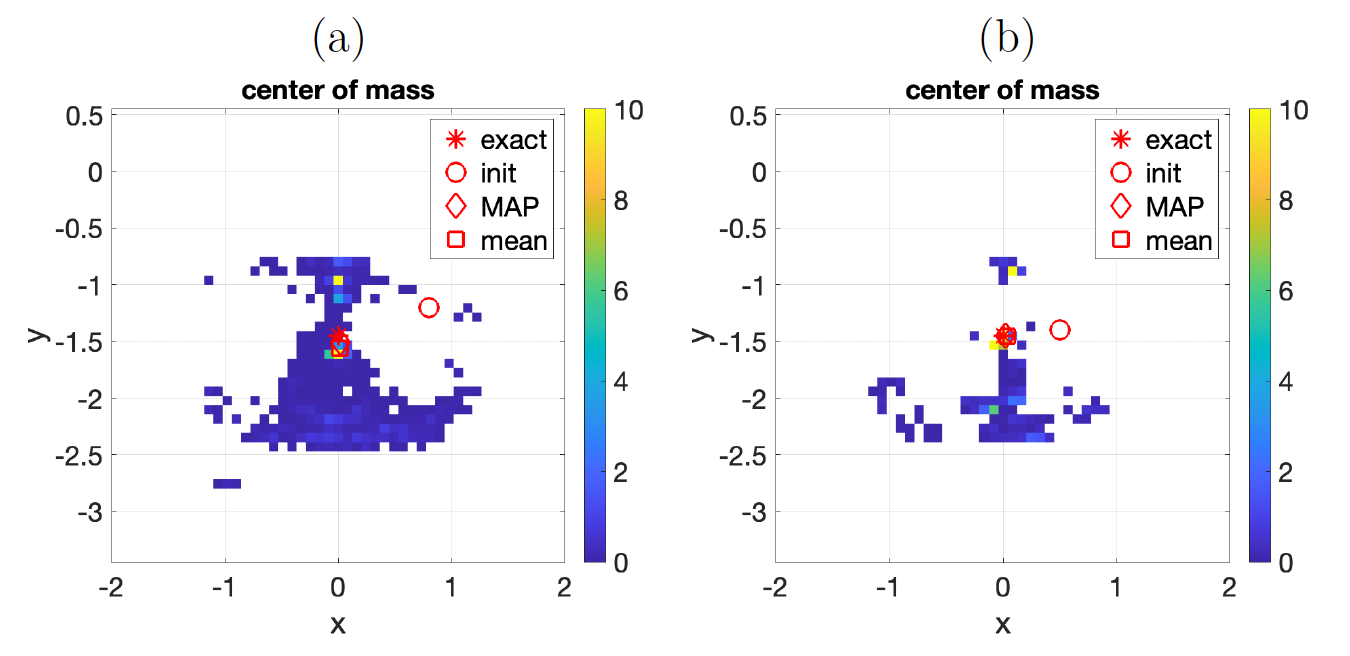} 
\caption{(a) Contour plot of the probability of belonging to the inclusion 
and (b) projection of the histogram for its material properties 
(pdf normalization) calculated by MCMC sampling and 
corresponding to Figure \ref{fig8}(a) with $5\%$ noise.}
\label{fig13}
\end{figure}

Figure \ref{fig8}(a)-(b) display sets of samples of the posterior distribution 
(\ref{unnormalized_ap}) with observation operator (\ref{observation_ap})
defined on a fixed uniform mesh for the data used in Figure \ref{fig4}.
The MAP estimates for the object geometry and its material properties are 
now provided by the sample with the highest probability. They are similar to 
the estimates obtained by optimization with either adaptive or fixed meshes
in Section \ref{sec:optimization}. We identify several families of samples in 
this figure. One wraps around the MAP estimate and the true object.
Other are similar to the local minima obtained with the algorithm in Section
\ref{sec:ofixed} on fixed meshes. The persistence of families of samples with 
different shape and material  properties, as well as the fact that their mean
departs from the MAP estimate suggest a multimodal posterior distribution. 
Increasing the signal to noise ratio by modifying the  properties of the bottom
layer we observe that multimodality features diminish, see panel (c). 

Panel (d) in Figure \ref{fig8} corresponds to a different inclusion and initial 
configuration. In this case, the dominant family of samples does not wrap around 
the MAP and the true object. This suggests that the prior information considered
is not enough, neither to convexify the cost, nor to properly approach the true 
object. This provides some insight on the optimization results in Figure \ref{fig7}: 
the adaptive scheme evolves towards the dominant mode, closer to the mean, 
because it has a larger basin of attraction. The scheme in Section \ref{sec:ofixed} 
is able to locate the MAP estimate provided by MCMC sampling. Notice that both 
employ fixed meshes and consider the same discretized observation operators 
and cost, whereas the adaptive optimization scheme considers a different 
approximation operator, thus a different cost functional.
Figure \ref{fig11}(a) represents the contour plot of the posterior
probability of a point belonging to this object. Panel 12(b)  and panel (a) in
Figure \ref{fig12} are projections of histograms representing joint material 
properties and the center of the inclusion, respectively. Two mains narrow 
peaks are identified in both, suggesting a bimodal distribution. 

\begin{figure}[!hbt]
\centering
\centering
\includegraphics[width=12cm]{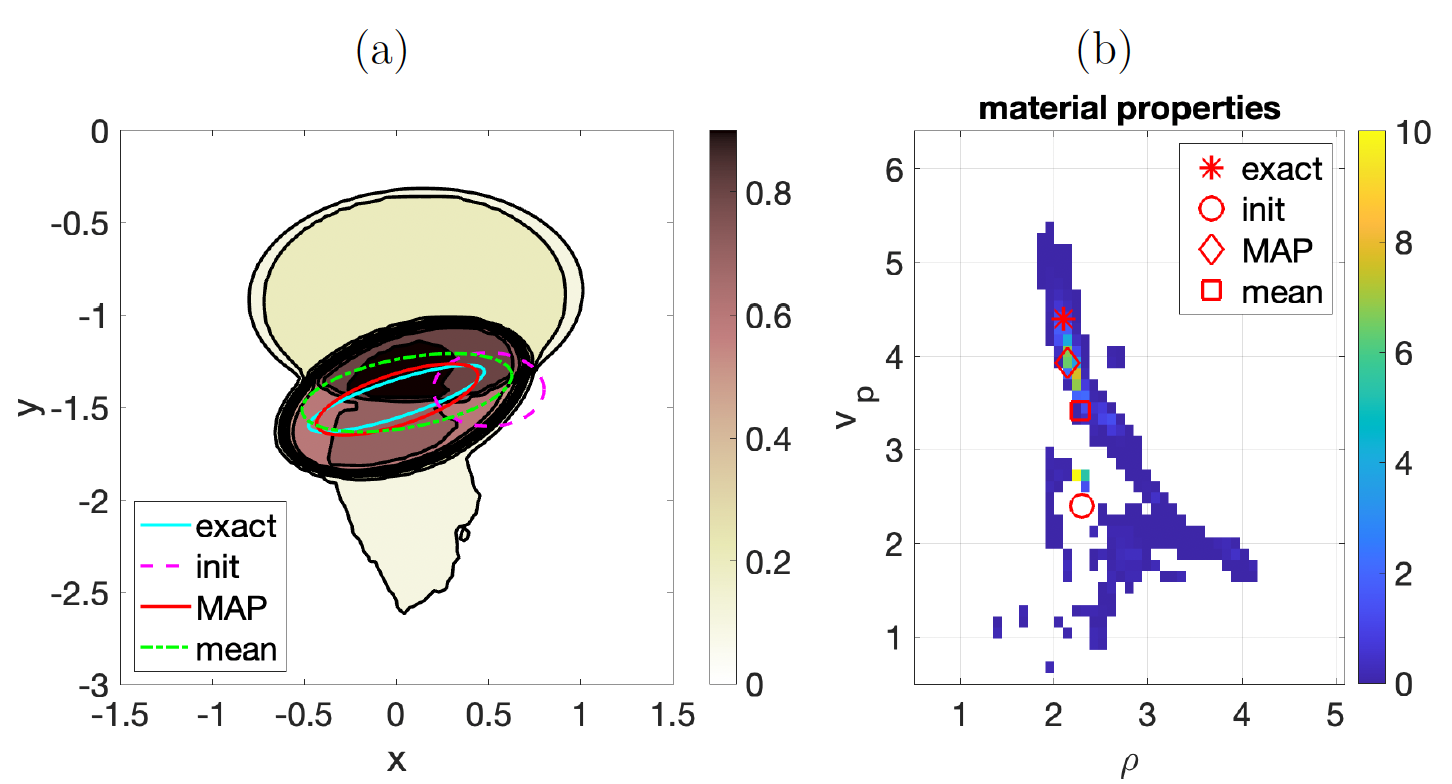} 
\caption{(a) Contour plot of the probability of belonging to the inclusion 
and (b) projection of the histogram for its material properties 
(pdf normalization) calculated by MCMC sampling and 
corresponding to Figure \ref{fig8}(b) with $15\%$ noise.}
\label{fig14}
\end{figure}

\begin{figure}[!hbt]
\centering
\includegraphics[width=12cm]{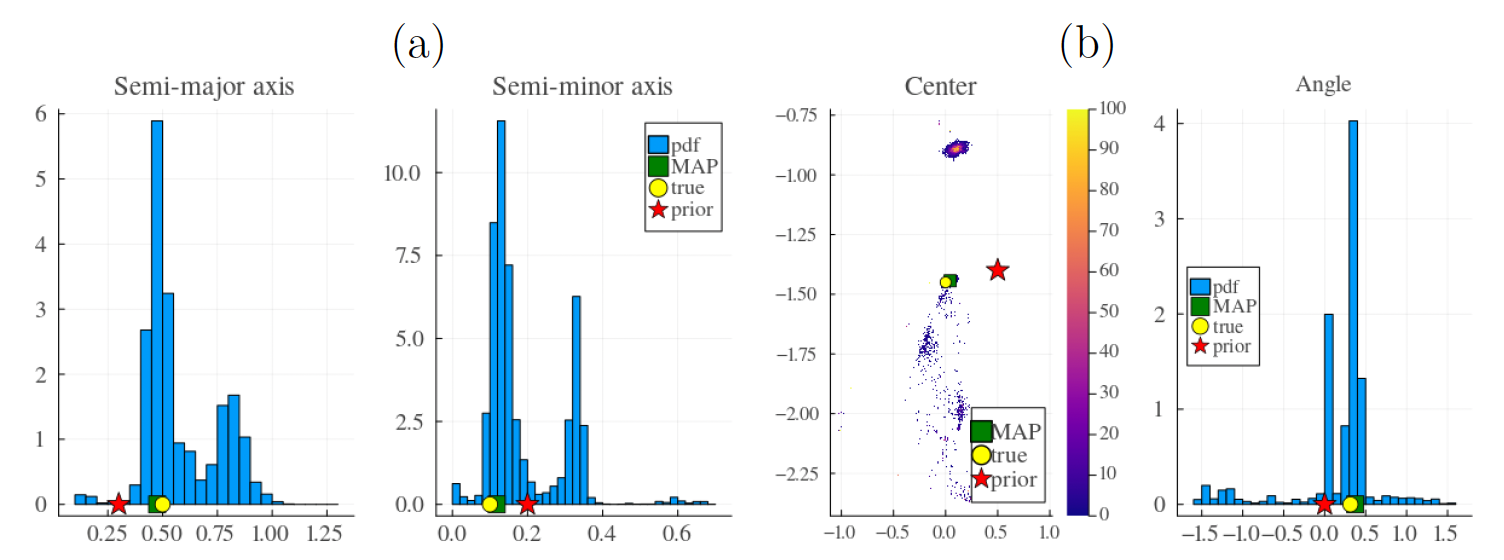} 
\caption{Histograms for axis size (a), centers and angles (b) associated to
Figure \ref{fig14} built from 412.800 samples.}
\label{fig15}
\end{figure}

Figure \ref{fig13} revisits the geometry in Figure \ref{fig1} with the $5\%$ noise
data considered in Section \ref{sec:optimization}. Panel (b) for the material properties
and panel (b) in Figure \ref{fig12} for the centers suggest the presence of three
narrow peaks, one of which is a dominant mode. 
Increasing the noise to $15\%$, Figures \ref{fig14} and \ref{fig15} suggest the
presence of two modes, visualizing the structure about them, the MAP estimates
and the mean values, compared to the true inclusion and the prior information.
One of the modes for the material properties seems reminiscent of the prior values 
while the other one relates to the true object that generated the data. 
Additional modes for $5\%$ may be related to the layered structure and the
difficulty of  resolving the depth.
These figures use uniform fixed meshes to evaluate the observation operator. 
Repeating the tests keeping the same data with stratified fixed meshes 
the main features remain similar. We may encounter small additional 
modes that vary with the mesh, specially for small noise levels in the data. 
However, the main modes remain similar.

We have repeated these simulations in larger computational domains, increasing
the number of receivers and sources, and also in homogeneous media. Working
with only one frequency, multiple solutions, ranked by their probability, persist.
Adding data obtained with additional frequencies, the main mode agrees for 
the different frequencies while the secondary modes vary. Resorting to costs
/likelihoods combining several frequencies, secondary modes are diluted. 


\section{Conclusions}
\label{sec:conclusions}

Devising computational tools to investigate subsurface structures from
surface measurements is a challenging problem, which we have addressed 
here using complementary strategies to infer details on the observed inclusions 
from synthetic data. Inclusions are defined by their material parameters and 
a few additional parameters describing their location, size and orientation. 
We have applied our algorithms for increasing noise levels in a simplified 
geometry of current interest involving localized salt inclusions. Salt areas 
typically hide reservoirs of relevant raw materials. 

In a Bayesian context we seek high probability configurations (MAP estimates)
by optimizing a cost regularized with prior information and constrained
with time dependent wave equations governing wave propagation.
First, we propose an automatic adaptive iterative scheme that employs  a
Levenberg-Marquardt-Fletcher type algorithm with observation operators 
defined by finite element solvers constructed on meshes that adapt 
automatically to the proposed objects at each iteration. Variations are
estimated using automatic differentiation. In synthetic tests employing a
single frequency, this  scheme evolves to the global minimum, which 
represents the MAP estimate of the sought inclusions, except in exceptional 
situations where additional minima with large attraction basins are encountered. 
We obtain basic uncertainty estimates in the parameter ranges by means 
of Laplace approximations. 
Second, we propose an alternative scheme that approximates observation 
operators resorting rougher solvers defined on either uniform or stratified 
meshes. These meshes are kept fixed during the process. We seek descent
directions estimating the variations of the solution as the inclusion changes 
by finite differences with small fixed steps. Depending on the selected step,
the algorithm may evolve to the global minimum or secondary local minima, if 
present.  

On one side, adaptive optimization employs more accurate solvers. On
the other, changing the mesh at each optimization step we may avoid getting 
trapped at spurious minima that might appear due to the usage of fixed meshes.
In practice, in the tests presented here, the configurations at which the
cost attains a global minimum remain similar, no matter whether we use uniform 
meshes (blind to the layered structure) or stratified meshes (adapted to the layers 
but not to the inclusion), for fine enough meshes. 
 
Markov Chain Monte Carlo methods have the potential of exploring the 
whole set of minima of the cost and the multimodal structure of the posterior 
probability, at the expense of solving a few million boundary value wave problems. 
Adapting the mesh to each proposed inclusion is not affordable. We are forced
to approximate the observation operator keeping a fixed mesh for different
inclusions.  We find MAP estimates similar to the ones provided by the adaptive 
optimization scheme in the tests we have performed, except in singular cases
in which the prior information allows for secondary minima with a larger
attraction basin.
This may be the case for elongated horizontal configurations, in which the lower 
half is screened by the upper half and the dataset we consider may not cover 
a wide enough range. 
MCMC studies provide a deeper understanding of uncertainty, since we identify 
additional high probability configurations. In the tests we have 
performed,  we have usually encountered  one mode related to the true inclusion,
and often a second mode that seems to keep a memory of the material parameters 
of the prior. 
A third intermediate mode appears occasionally, using either uniform or stratified 
meshes. Small fluctuations may vary with the mesh, specially when the 
noise level in the data is small.
However, the dominant modes remain similar for the two types of meshes, uniform
and stratified, provided they are fine enough, though the latter is expected to reproduce
better transmission of waves at interfaces. This is a useful remark since, in practice,
we may not know the underlying layered structure and be forced to use blind
uniform meshes. This is often the case in high-dimensional studies.
This situation  remains even if we enlarge the computational region and increase the 
number of  emitters/receivers or we suppress the layers and consider a uniform 
background.

Notice that our tests use synthetic data generated numerically. Therefore we know 
beforehand the true inclusions that produce the data. In the tests we have performed 
the MAP estimate  is related to the true object. However, in the absence of a global 
view of the probability modes that may appear, there is the risk of trusting secondary
modes as inaccurate reconstructions. Global uncertainty studies are thus important 
to clarify these facts and may be useful to guide the design of imaging set-ups 
for which multimodality is suppressed and uncertainty in the predictions of inclusion 
properties diminishes.  In our set-up,
varying the frequency of the emitted signals, the secondary modes vary while the
main mode remains. Thus, considering costs/likelihoods that incorporate data 
recorded with different frequencies,  we may suppress secondary minima/modes.
Then the results provided by MCMC studies would be fully similar to those 
automatically obtained by the optimization algorithms and Laplace approximations.

We have focused here on low dimensional descriptions of the inclusions, in the
sense that their geometry is described by a few fixed parameters. More precise
descriptions involve coefficients belonging to infinite dimensional spaces, often
approximated for computational purposes by their values on meshes. This leads
to high dimensional problems that can also be treated by optimization and MCMC
techniques. Conclusions extracted from low dimensional studies about set-up
design to reduce uncertainty and meshes can be useful when dealing with
high dimensional formulations.


\vskip 2mm
{\bf Acknowledgements.}
This research has been partially supported by Spanish FEDER/MICINN - 
AEI grant No. PID2020-112796RB-C21 and by a  Doctoral Scholarship Agencia 
Nacional de Investigaci\'on y Desarrollo, 2022, Chile (CA). 



\section*{Appendix A: Nondimensionalization and parameters}
\label{sec:nondimensional}

To write the forward problem in dimensionless variables we 
choose two characteristic time $T$ and length $L$ scales.
We set $\mathbf x= \tilde {\mathbf x} L$, $t= \tilde t T$,  $u = \tilde u L $, 
$R =  \tilde R L$ and $\Sigma = \tilde \Sigma L$. 
Let us denote
$\tilde \rho = \frac{\rho}{\rho_0}$, $\tilde v_p= \frac{v_p T}{L}$,
$\tilde f(\tilde t)= \frac{T^2}{L} f_0 (1- 2 \pi^2 (f_M T)^2 \tilde t^2) 
e^{-\pi^2 (f_M T)^2 \tilde t^2}$ and 
$\tilde G(\tilde {\mathbf x}) =  \frac{1}{ ( \pi \kappa )^{n/2} } \sum_{k=1}^K
\exp(- \frac{|\tilde {\mathbf x} - \mathbf x_k|^2}{\kappa })$, $n=2$.
We set  $T=1$ s, $L= 1000$ m  and $\rho_0= 1000$ kg/m$^3$.
We choose  $f_M = 2$, $f_0 = 100$ so that
$f_M T = 2$ and $f_0 \frac{ T^2 }{L} = 0.1$.
Making the change of variables and dropping the $ \tilde{} $ for ease 
of notation, we get
\begin{eqnarray}
\begin{array}{ll}
\rho  u_{tt} -  {\rm div}(\rho v_p^2 \nabla u) =  \rho  f(t) 
G(\mathbf x),  & \mathbf x \in R, \; t>0, \\
\nabla u \cdot \mathbf n = 0,  & \mathbf x \in \Sigma, \; t>0, \\ 
\frac{\partial u}{\partial \mathbf n} = - \frac{1}{v_p} u_t, &
\mathbf x \in \partial R- \Sigma, \, t>0, \\
u(\mathbf x, 0) = 0, u_t(\mathbf x, 0) = 0, & \mathbf x \in R,
\end{array} \label{forward_adim2}
\end{eqnarray}
where $R$ is the computational rectangular region, $\Sigma$
the upper surface and $v_p$ stands for the local wave velocity.
On  $\partial R- \Sigma$ we enforce an approximate nonreflecting 
boundary condition.

Table \ref{table1} collects dimensionless parameter values for Figure
\ref{fig1}.
The dimensionless computational region becomes  $R=[-1.5, 1.5] 
\times [-3,0]$. We discretize the problem using  FEM meshes of  
minimum step $0.04$
with emitters/receivers interspaced with step $0.02$. Emitters are 
$e_k=-1+0.04k$, $k=0,\ldots,50$, and receivers
$r_j=-1.02+0.04j$, $j=0,\ldots,51$.  
The value $\kappa$ is adjusted to the mesh so that it affects a small
region about the emitters.  Here we have set $\kappa=0.04$. 
The time step for the numerical method is $\delta t = 10^{-3}$
and the profiles are recorded at intervals of $10^{-1}$. The maximum time
is $T= 2.5$.

\begin{table}[!hbt]
\centering
\begin{tabular}{|c|c|c|c|c|c|c|}
\hline
 {\rm Layers} & 1 & 2 & 3 & 4 & 5 & {\rm Object} \\
 \hline
 $\rho $ & 2 & 2.5 & 2.49 & 2.49 &  2.6  & 2.1   \\
 \hline
 $v_p $ & 1.5 & 2.5 & 2.8 & 3.3 &   3.1 & 4.4  \\
\hline
\end{tabular}
\caption{Dimensionless parameter values for the true layered geometry.}
\label{table1}
\end{table}

We generate synthetic data  $\mathbf d_{\rm true}$ for the numerical tests solving 
numerically (\ref{forward_adim2}) and evaluating the approximate solution at the 
detectors in a fixed time grid. Approximation schemes are described in Appendix
C. Then, we corrupt the data with noise to obtain  $\mathbf d$ according to
\begin{eqnarray*}
d_j^m = d_{j,\rm true}^m + {r \over 100 }\sigma \beta,
\; j=1,\ldots,J, \, m=1,\ldots,M, 
\; \; 
\sigma =  \Big(\sum_{j=1}^J \sum_{m=1}^M 
{ |d_{j,\rm true}^m|^2  \over JM } \Big)^{1\over 2}
\end{eqnarray*}
where $\beta$ is drawn from ${\cal N}(0,1)$ and $r>0$ is the
noise level. Then $\sigma_{\rm noise} = \sigma \, r/100$.


\section{Appendix B: Well posedness of the truncated forward problem}
\label{sec:truncated_wp}


We establish next existence, uniqueness and regularity results for solutions of
(\ref{balance_escalar}) when $\gamma >0$.
We assume that $R \subset \mathbb R^n$, $n \geq 2$, is a truncated half-space 
with borders defined by hyperplanes, that is, a rectangle when $n=2$ or a
parallelepiped when $n=3$. 

We denote by $H^m(R)$ and $L^2(R)$ the standard Sobolev spaces and the space
of square-integrable functions, respectively. 
$L^2(\partial R \setminus \overline{\Sigma})$ is the space of traces on the boundary
\cite{lionsmagenes,grisvard}. 
Similarly, we denote by $C^m([0,T]; H)$, $m \geq 0$, and $L^\infty([0,T]; H)$ the
spaces of continuously differentiable functions up to the order $m$ and bounded
functions with values in a Hilbert space $H$, respectively \cite{lionsmagenes}.
For any $u(t) \in H^1(R)$ with $u_t(t) \in L^2(R)$, we define the energy as
\begin{eqnarray}
E(u(t),\!u_t(t)) = {1\over 2} \int_{R} \ \rho |u_t|^2  d \mathbf x  + 
{1\over 2} \int_{R}  \chi |\nabla u|^2   d \mathbf x
+ {1\over 2} \! \int_{\partial R \setminus \overline{\Sigma}} 
\mu \chi \gamma |u|^2   d S_{\mathbf x}.
\label{en_def}
\end{eqnarray} 

Problem (\ref{balance_escalar}) admits a weak formulation.
Formally, multiplying by $w \in H^1(R)$, integrating by parts over $[0,T] \times R$ 
and assuming that $u$ is smooth enough, we find
\begin{eqnarray} \begin{array}{r}
\displaystyle
{d \over dt^2}\int_{R}    \rho(\mathbf x)   u(t,\mathbf x)  w(\mathbf x) \, d \mathbf x 
+ \int_{R}  \chi(\mathbf x) \nabla u(t,\mathbf x)  \nabla w(\mathbf x) \, d \mathbf x   +
\\ [2ex] \displaystyle
{d \over dt}\int_{\partial R \setminus \overline{\Sigma}}  \chi(\mathbf x) \gamma(\mathbf x) 
u(t,\mathbf x) w(\mathbf x) \, d S_{\mathbf x}   
 =     \int_{R}   \rho(\mathbf x) h(t,\mathbf x) w(\mathbf x)  \, d \mathbf x
\\[2ex] \displaystyle 
u(0) =  u_0, \quad u_t(0) = u_1, 
\end{array} \label{balance_weak}
\end{eqnarray}
for all $w \in H^1(R)$, given $f \in L^\infty(0,T;L^2(R))$.
We seek a solution $u$ to (\ref{balance_weak}), at least with regularity 
$C([0,T]; H^1(R)) \cap C^1([0,T]; L^2(R))$ 
to recover (\ref{balance_escalar}) in the sense of distributions.

{\bf Theorem 1.} {\it Let us assume that 
\begin{itemize}
\item $\rho, \chi, \alpha \in L^\infty(R)$, $0 < \rho_{\rm min} \leq \rho \leq \rho_{\rm max}$, 
$0 < \chi_{\rm min} \leq \chi \leq \chi_{\rm max}$, $0 < \gamma_{\rm min} \leq 
\gamma \leq \gamma_{\rm max}$, 
\item  $u_0 \in H^1(R)$, $u_1 \in L^2(R)$, $h \in C([0,T];L^2(\Omega))$.
\end{itemize}
Then, there exists a unique solution 
$u \in C([0,T]; H^1(R)) \cap C^1([0,T]; L^2(R)) $ to problem (\ref{balance_weak})  
with
$u_t \in L^2(0,T;L^2(\partial R \setminus \overline{\Sigma}))$. 
This solution satisfies the wave equation in the sense of distributions.
Moreover, it satisfies an energy inequality and depends continuously on the data. 
The following estimates hold
\begin{eqnarray} \begin{array}{l} \displaystyle
 \|u_t\|_{L^\infty(0,T;L^2(R))} \leq K(T, \rho_{\rm min}, \rho_{\rm \max}, E(u_0,u_1), 
 \|h\|_{C([0,T];L^2(R))}),  \\[1ex]
 \| \nabla u \|_{L^\infty(0,T;L^2(R))} \leq K(\chi_{\rm min}, T, \rho_{\rm min}, \rho_{\rm \max}, 
 E(u_0,u_1), \|h\|_{C([0,T];L^2(R))}),  \\[1ex]
 \|u_t\|_{L^2(0,T;L^2(\partial R \setminus \overline{\Sigma}))} \leq 
 K(\gamma_{\rm min}, \chi_{\rm min}, T, \rho_{\rm min}, \rho_{\rm \max}, E(u_0,u_1), 
 \|h\|_{C([0,T];L^2(R))}),  \\[1ex]
 \|u\|_{L^\infty(0,T;L^2(R))} \leq K(\mu,  T, \rho_{\rm min}, \rho_{\rm \max}, E(u_0,u_1), 
 \|h\|_{C([0,T];L^2(R))}),  
\end{array} \label{continuous1} \end{eqnarray}
for any $\mu>0$, where $K(\cdot)$ denote different positive constants depending 
continuously on the specified arguments.}

{\bf Proof.} The proof is based on the use of Galerkin bases and compactness arguments.
When $\gamma=0$, we can consider a Galerkin basis formed by eigenfunctions
of an elliptic operator and prove explicit strong convergence results for the
eigenfunction expansion of the solution \cite{raviart}. However, if $\gamma>0$,
this approach fails, thus, we resort to general abstract bases and compactness 
arguments next.

{\it Step 1: Galerkin approximation.} Since $H^1(R)$ is separable we can always
find a set $\{\phi_1, \ldots, \phi_k, \ldots \}  \subset H^1(R)$ whose elements
are linearly independent  while their linear combinations with real coefficients are dense 
in $H^1(R)$ \cite{lionsmagenes}.
For each $M \in \mathbb N$, we denote by $V^M$ the space generated by 
$\{\phi_1, \phi_2, \ldots, \phi_M\}$ and consider the approximate problem: Find
$u^M(t, \mathbf x) = \sum_{m=1}^M a_m(t) \phi_m(\mathbf x)$ 
such that
\begin{eqnarray} 
\begin{array}{l} 
{d^2\over dt^2} \int_{R}  \rho u^M(t)  w  \, d \mathbf x +
\int_{R}  \chi \nabla  u^M(t)  \nabla  w  \, d \mathbf x   
+ {d \over dt} \int_{\partial R \setminus \overline{\Sigma}} 
\chi \gamma \, u^M(t)  w \, d S_{\mathbf x}  \\
 \displaystyle  \hfill =     \int_{R} \rho h(t)  w \,  d \mathbf x, 
\\[2ex] 
u^M(0) =  u_0^M, \quad u_t^M(0) =  u_1^M. 
\end{array}  \label{eq:varwM}
\end{eqnarray}
for all $w \in V^M$  and $t \in [0,T]$, see \cite{lionsmagenes}, where
$u_0^M= \sum_{m=1}^M u_{0,m} \phi_m$ and 
$u_1^M= \sum_{m=1}^M u_{1,m} \phi_m$  
are the projections of $u_0$ and $u_1$ in $V^M$.

{\it Step 2: Change of variables.}
To achieve the necessary estimates, we change variables and set
$u^M = e^{\mu t} v^M$, $\mu >0$, so that
$u^M_t = \mu e^{\mu t} v^M + e^{\mu t} v^M_t$ and 
$u^M_{tt} = \mu^2 e^{\mu t} v^M  + 2 \mu e^{\mu t} v^M_t + e^{\mu t} v^M_{tt}$. 
Problem (\ref{eq:varwM}) becomes: Find $v^M =
\sum_{m=1}^M b_m(t) \phi_m(\mathbf x)$ such that
\begin{eqnarray} 
\begin{array}{l} 
{d^2\over dt^2} \int_{R}  \rho v^M(t)  w  \, d \mathbf x +
\int_{R}  \chi \nabla  v^M(t)  \nabla  w  \, d \mathbf x   
+ {d \over dt} \int_{\partial R \setminus \overline{\Sigma}} 
\chi \gamma \, v^M(t)  w \, d S_{\mathbf x} 
\\[2ex] 
+  \int_{R}  \rho \mu^2  v^M(t)  w  \, d \mathbf x
+ {d \over dt} \int_{R}  2 \rho  \mu  v^M(t)  w  \, d \mathbf x
+ \int_{\partial R \setminus \overline{\Sigma}} \chi \gamma
 \mu v^M(t)  w \, d S_{\mathbf x} 
\\[2ex] 
 \hfill = e^{-\mu t} \ \int_{\Omega}  \rho h(t)   w \,  d \mathbf x, 
\\[2ex] 
v^M(0) =  u_0^M, \quad v_t^M(0) =  u_1^M, 
\end{array}  \label{eq:varwM_lam}
\end{eqnarray}
for all $w \in V^M$ and $t \in [0,T]$.

{\it Step 3: Existence of an approximant.} Problem (\ref{eq:varwM_lam}) is
equivalent to a linear system of $M$ second order differential equations for 
the coefficient functions $b_m$
\begin{eqnarray}  
\begin{array}{l} 
\hfill \sum_{m=1}^M b_m''(t) \int_R \rho \phi_m \phi_k \, d \mathbf x 
\!+\! \sum_{m=1}^M b_m'(t)  \left[ 2 \mu \int_{R}  \rho \phi_m  \phi_k   \, d \mathbf x 
\!+\! \int_{\partial R \setminus \overline{\Sigma}} \chi \gamma \,\phi_m  \phi_k \, d S_{\mathbf x}   
\right] 
\\[2ex] 
+ \sum_{m=1}^M   b_m(t)  \left[  \int_{R} \chi \nabla \phi_m \nabla \phi_k   \, d \mathbf x    
+  \mu^2 \int_{R}  \rho  \phi_m  \phi_k    \, d \mathbf x +  \mu \int_{\partial R \setminus 
\overline{\Sigma}} \chi \gamma \phi_m  \phi_k   \, d S_{\mathbf x} 
\right] 
\\[2ex] 
\hfill =  e^{-\mu t}    \int_{\Omega}  \rho h(t) \phi_k \,  d \mathbf x,   
\\[2ex] 
b_m(0) =   u_{0,m}, \quad
b_m'(0) =  u_{1,m},  \quad m=1, \ldots, M,  
\end{array}  \label{eq:varwM_coef} 
\end{eqnarray}
for $k=1,\ldots, M$. In matricial form,
\begin{eqnarray*}
\mathbf M \mathbf b'' + \mathbf D \mathbf b' + \mathbf A \mathbf b = \mathbf h(t),
\end{eqnarray*}
where $\mathbf h(t) \in C([0,T])$. This linear system can be written as a first
order linear system for $\mathbf b$ and $\mathbf a = \mathbf b'$, which has a 
unique classical solution $\mathbf b = (b_1,\ldots, b_m) \in C^2([0,T])$
for any $M$, see \cite{coddington}, Ch. 3.3.

{\it Step 4: Uniform estimates.}
We multiply (\ref{eq:varwM_coef}) by $b_k'$  and add over $k$ to get
\begin{eqnarray} 
\begin{array}{l} 
{1\over 2}{d\over dt} \int_{R} \rho |v^M_t(t)|^2  d \mathbf x  
+ \int_{\partial R \setminus \overline{\Sigma}}  \chi \gamma |v^M_t(t)|^2  d S_{\mathbf x} 
+ 2 \mu \int_{\Omega} \rho |v^M_t(t)|^2  d \mathbf x  +
\\[2ex]  
{1\over 2}{d\over dt}\left [ \int_{R} \chi |\nabla v^M(t)|^2  d \mathbf x  
+  \mu^2 \int_{R}  \rho  |v^M(t)|^2    \, d \mathbf x
+  \mu \int_{\partial R \setminus \overline{\Sigma}} \chi \gamma |v^M(t)|^2  \, d S_{\mathbf x} 
\right]
\\[2ex] 
=   e^{-\mu t}  \int_{\Omega} \rho h(t) v^M_t(t)  \,  d \mathbf x.  
\end{array}  \label{eq:denergyM}
\end{eqnarray}
For any $v(t) \in H^1(R)$ with $v_t(t) \in L^2(R)$, we define the energy as
\begin{eqnarray} \begin{array}{l}
E_\mu(u(t),\!u_t(t)) \!=\! {1\over 2} \! \int_{R}   \rho |u_t|^2  d \mathbf x  \!+\!  {1\over 2} 
\! \int_{R}    \left[ \chi |\nabla u|^2  \!+\!  \mu^2   \rho  |u|^2  \right]   d \mathbf x
\!+\! {1\over 2} \! \int_{\partial R \setminus \overline{\Sigma}} 
  \mu \chi \gamma |u|^2   d S_{\mathbf x}.
\end{array} \label{en_def_mu}
\end{eqnarray} 
Integrating (\ref{eq:denergyM}) over $t$ it follows that
\begin{eqnarray} 
\begin{array}{l}
 E_\mu(v^M(t),v^M_t(t))  + \int_0^t \int_{\partial R \setminus \overline{\Sigma}}  
 \chi \gamma |v^M_t(s)|^2  d \mathbf x ds
 + 2\mu \int_0^t \int_{R} \rho |v^M_t(s)|^2  d \mathbf x ds
\\[2ex]  
= E_\mu(u_0^M,u_1^M) + \int_0^t \int_{R} e^{-\mu s} \rho  h(s)  v^M_t(s)  \,  d \mathbf x ds.
\end{array}  \label{eq:ienergyM}
\end{eqnarray}
Discarding positive terms and using the properties of $\rho$, we find
\begin{eqnarray} 
\begin{array}{ll} 
\rho_{\rm min}  \|v_t^M(t)\|_{L^2(R)}^2 \leq &  2 E_\mu(u_0^M,u_1^M) + 
 \rho_{\rm max}   \int_0^t  \|h(s)\|_{L^2(R)}^2 ds
\\[2ex]  
& +  \rho_{\rm max}  \int_0^t \|v_t^M(s)\|_{L^2(R)}^2 ds, 
\end{array}  \label{eq:ivtM}
\end{eqnarray}
thanks to Young's inequality. Notice that $E_\mu(u_0^M,u_1^M)
\rightarrow E_\mu(u_0,u_1)$ as $M \rightarrow \infty$ due to strong
convergence in $H^1(R)$ and $L^2(R)$.
Then, Gronwalls' inequality \cite{haraux}  yields a uniform bound on
$\|v_t^M\|_{L^\infty(0,T;L^2(R))}$ in terms of $T$, $\|h\|$, $E(u_0,u_1)$, 
and $\rho$. Inserting this uniform estimate in inequality (\ref{eq:ienergyM})
we obtain uniform bounds on 
$ \|v_t^M\|_{L^2(0,T;L^2(\partial R \setminus \overline{\Sigma}))}$,
$ \|v^M\|_{L^\infty(0,T;L^2(\partial R \setminus \overline{\Sigma}))}$
and 
$ \|v^M\|_{L^\infty(0,T;H^1(R))}$ when $\mu>0$.

{\it Step 5: Compactness.}
By classical compactness results \cite{lionsmagenes,grisvard}, we can extract a 
subsequence $v^{M'}$ converging weakly star in $W^{1,\infty}(0,T;L^2(R)) 
\cap L^\infty(0,T;H^1(R))$ to a limit 
\[ v \in W^{1,\infty}(0,T;L^2(R)) \cap L^\infty(0,T;H^1(R)) \]
as $M' \rightarrow \infty$, with traces  
$v^{M'}_t\big|_{\partial R \setminus \overline{\Sigma}}$ converging weakly in 
$L^2(0,T;L^2(\partial R \setminus \overline{\Sigma}))$ to a limit
$v_t\big|_{\partial R \setminus \overline{\Sigma}}$ and 
$v^{M'}\big|_{\partial R \setminus \overline{\Sigma}}$ converging weakly star
$L^\infty(0,T;L^2(\partial R \setminus \overline{\Sigma}))$ to a limit
$v\big|_{\partial R \setminus \overline{\Sigma}}$.
Moreover, ${d^2 \over dt^2 }v^{M'} \rightarrow {d^2 \over dt^2} v $ in the
sense of distributions.


{\it Step 6: Passage to the limit.} 
To find the equation satisfied by $u$, we take $w=\phi_k$, multiply 
(\ref{eq:varwM_lam}) by $\psi(t) \in C_c^\infty([0,T))$ and integrate over $t$  
to obtain
\begin{eqnarray*} 
\begin{array}{r} 
\int_0^T \int_{R} \rho  v^{M'}   \psi_{tt} \phi_k \, d \mathbf x 
ds +  \int_{R} \! \rho  u_{1,m}  \psi(0) \phi_k \, d \mathbf x 
-  \int_{R} \! \rho  u_{0,m}   \psi_t(0)  \phi_k \, d \mathbf x 
\\[2ex] 
+ \int_0^T \int_{R}  \chi \nabla  v^{M'}  \nabla  \phi_k \psi  \, d \mathbf x  ds 
+ \int_0^T \int_{\partial R \setminus \overline{\Sigma}} \chi \gamma v^{M'}_t  
\phi_k \psi \, d S_{\mathbf x}  ds
\\[2ex] 
+ \int_0^T \int_{R}  \rho \mu^2  v^{M'}  \phi_k \psi  \, d \mathbf x ds
+ \int_0^T \int_{R}  2 \rho  \mu  v^{M'}_t  \phi_k \psi   \, d \mathbf x ds
\\[2ex] 
+ \int_0^T \int_{\partial R \setminus \overline{\Sigma}} 
\chi \gamma \mu v^{M'}  \phi_k \psi  \, d S_{\mathbf x}  ds
= \int_0^T e^{-\mu s}   \int_{R}  h(s) \phi_k \psi  \,  d \mathbf x ds, 
\end{array}  
\end{eqnarray*}
for $k \leq M'$. 
Letting $M' \rightarrow \infty$ we find  
\begin{eqnarray} 
\begin{array}{r} 
\int_0^T  \int_{R} \! v    \psi_{tt} \phi_k \, d \mathbf x 
ds +  \int_{R} \! \rho  u_1  \psi(0) \phi_k \, d \mathbf x 
-  \int_{R} \! \rho  u_0  \psi_t(0)  \phi_k \, d \mathbf x 
\\[2ex] 
+ \int_0^T \int_{R}  \chi \nabla  v   \nabla  \phi_k \psi  \, d \mathbf x  
ds
+ \int_0^T \int_{\partial R \setminus \overline{\Sigma}} \chi \gamma 
v_t  \phi_k \psi \, d S_{\mathbf x}  ds
\\[2ex] 
+ \int_0^T \int_{R}  \rho \mu^2  v   \phi_k \psi  \, d \mathbf x ds
+ \int_0^T \int_{R}  2 \rho  \mu  v_t  \phi_k \psi   \, d \mathbf x ds
\\[2ex] 
+ \int_0^T \int_{\partial R \setminus \overline{\Sigma}} \chi \gamma
\mu v   \phi_k \psi  \, d S_{\mathbf x}  ds
= \int_0^T e^{-\mu s} \int_{R}  h(s)  \phi_k \psi  \,  d \mathbf x ds, 
\end{array}  
\label{eq:v}
\end{eqnarray}
for all $\phi_k$. The identity extends to all $w \in H^1(R)$ by density. 
Taking $\psi \in C_c(0,T)$ and $\phi \in C_c(R)$ in (\ref{eq:v}), and integrating
by parts, we see that $v$ satisfies the equation
$\rho v_{tt} -  {\rm div}(\chi \nabla v) + 2 \rho \mu v_t + \rho \mu^2  v
=  e^{-\mu t} h$ in the sense of distributions in $[0,T] \times R$, see 
\cite{lionsmagenes}.
Undoing the change of variables, we have constructed a solution $u$ of
\begin{eqnarray}
\rho u_{tt} -  {\rm div}(\chi \nabla u) =  h \quad {\rm in} \,  {\cal D}'(0,T) \times R
\label{eq:u}
\end{eqnarray}
in the sense of distributions.  

Since $u \in L^2(0,T;H^1(R))$, $u_t \in L^2(0,T;L^2(R))$ and 
$u_{tt} \in L^2(0,T;(H(R))')$, after eventually modifying a set of zero 
measure,  $u \in C([0,T];H^1(R))$ and $u_t \in C([0,T];L^2(R))$, 
see Theorem 8.2 in \cite{lionsmagenes}, Ch. 3. Then, 
$u(0) \in H^1(R)$ and $u_t(0) \in L^2(R)$.
We take $\psi \in C([0,T))$ and $\phi \in C_c(R)$ in (\ref{eq:v}), 
integrate by parts, and use (\ref{eq:u}), to get $u(0)=u_0$ and
$u_t(0)=u_1$. Therefore, we have constructed a solution 
$u \in C([0,T]; H^1(R)) \cap C^1([0,T]; L^2(R))$ 
to (\ref{balance_escalar}) with
$u_t \big|_{\partial R \setminus \overline{\Sigma}} \in 
L^2(0,T; L^2(\partial R \setminus \overline{\Sigma})).$

{\it Step 7: Energy inequality.} Taking limits in (\ref{eq:ienergyM}) and using
the properties of weak convergences we find
\begin{eqnarray} 
\begin{array}{l} 
 E_\mu(v(t),v_t(t))  + \int_0^t \int_{\partial R \setminus \overline{\Sigma}}  
 \chi \gamma |v_t(s)|^2  d \mathbf x ds
 + 2\mu \int_0^t \int_{R} \rho |v_t(s)|^2  d \mathbf x ds
\\[2ex]  
\leq E_\mu(u_0,u_1) + \int_0^t \int_{R} e^{-\mu s} \rho  h(s)  v_t(s)  \,  d \mathbf x ds.
\end{array}  \label{eq:ienergyv}
\end{eqnarray}
Undoing the initial change of variables, the function $u$ satisfies a similar 
inequality with $\mu=0$:
\begin{eqnarray} 
\begin{array}{l} 
\int_{R}  \rho |u_t(t)|^2  d \mathbf x  + \int_{R} \chi |\nabla u(t)|^2  \, d \mathbf x
+ 2 \int_0^t \int_{\partial R \setminus \overline{\Sigma}}  \chi \gamma
|u_t(s)|^2  d \mathbf x ds
\\[2ex]  
\leq  \int_{R}  \rho |u_1|^2  d \mathbf x  + \int_{R} \chi |\nabla u_0|^2   
\, d \mathbf x + 2 \int_0^t \int_{R}  h(s) g u_t(s)  \,  d \mathbf x ds.
\end{array}  \label{eq:ienergyu}
\end{eqnarray}

{\it Step 8: Dependence on parameters.}
Inequality (\ref{eq:ienergyu}) implies
\begin{eqnarray} 
\begin{array}{ll} 
\rho_{\rm min}  \|u_t(t)\|_{L^2(R)}^2 \leq &  2 E(u_0,u_1) + 
 \rho_{\rm max}   \int_0^t  \|h(s)\|_{L^2(R)}^2 ds
\\[2ex]  
& +  \rho_{\rm max}  \int_0^t \|u_t(s)\|_{L^2(R)}^2 ds. 
\end{array}  \label{eq:iutM}
\end{eqnarray}
By Gronwall's inequality,  
$\|u_t\|_{L^\infty(0,T;L^2(R))}$ is bounded in terms of $\rho_{\rm min}$, 
$\rho_{\rm max}$, $T$, $E(u_0,u_1)$, and $\|f\|_{L^\infty(0,T;L^2(R))}$.
Inserting this information in (\ref{eq:ienergyu}) we obtain similar
estimates for $\|\nabla u\|_{L^\infty(0,T;L^2(R))}$ and 
$\|u_t\|_{L^2(0,T;L^2(\partial R \setminus \overline{\Sigma}))}$ with 
constants depending also on $\chi_{\rm min}$ and $\gamma_{\rm min}$. 
If we wish to estimate
$\|u\|_{L^\infty(0,T;L^2(R))}$ we need to take limits in 
(\ref{eq:ienergyM}) and (\ref{eq:ivtM}) with $\mu>0$ to bound
$\|v\|_{L^\infty(0,T;L^2(R))}$, $v= e^{-\mu t} u$, and then
$\|u\|_{L^(0,T;L^2(R))}$. The bounding constant depends
now also on $\mu$.

{\it Step 9: Uniqueness.} Let us consider two solutions $w_1$ and $w_2$
with the stated regularity. Then, $u=w_1-w_2$ is a solution of a similar problem
with zero source and zero initial conditions. We perform the change of variables
$v= e^{-\mu t} u$ and consider the equation (\ref{eq:v}) for $v$. 
Next, we adapt the uniqueness  proof of Theorem 8.1 in \cite{lionsmagenes}, 
Ch. 3. We choose $\tau \in (0,T)$ and take a test function $w(t)  = - \int_t^\tau 
v(\sigma) d \sigma$ when $t <\tau$, zero otherwise. With this choice, 
$w_t = v \in C([0,T],H^1(R))$ and $w \in C([0,T],H^1(R))$.
Inserting the test function $w$ in (\ref{eq:v}) instead of $\phi \psi$ we get
\begin{eqnarray*} 
\begin{array}{r}   
\int_0^T\int_R \left[- \rho v_{t} w_t - 2 \rho \mu v w_t + \chi \nabla v
\nabla w + \rho \mu^2  v w \right] d \mathbf x ds   \\
[2ex]  
+ \int_0^T\int_{\partial R \setminus \overline{\Sigma}}
[- \chi \gamma v  w_t + \chi \gamma \mu v  w] d S_{\mathbf x} ds  = 0,
\end{array}  
\end{eqnarray*}
As a result, 
\begin{eqnarray*} 
\begin{array}{r}   
\int_0^\tau \int_R \left[- \rho v_{t} v - 2 \rho \mu |v|^2 + \chi \nabla w_t
\nabla w + \rho \mu^2  w_t w \right] d \mathbf x ds   \\
[2ex]  
+ \int_0^\tau \int_{\partial R \setminus \overline{\Sigma}}
[- \chi \gamma |v|^2 + \chi \gamma \mu w_t  w] d S_{\mathbf x} ds  = 0.
\end{array}  
\end{eqnarray*}
Integrating we find
\begin{eqnarray*} 
\begin{array}{r}  
0= {1\over 2} \int_R \left[- \rho |v(\tau)|^2 - \chi |\nabla w(0)|^2- \rho \mu^2  |w(0)|^2
- 2 \rho \mu |v|^2  \right] d \mathbf x ds   
\\ [2ex]  
- {1\over 2} \int_{\partial R \setminus \overline{\Sigma}}
\chi \gamma \mu w(0) d S_{\mathbf x} ds 
- \int_0^\tau \int_R 2 \rho \mu |v|^2 d \mathbf x ds  
- \int_0^\tau \int_{\partial R \setminus \overline{\Sigma}}
  \chi \gamma |v|^2  d S_{\mathbf x} ds.
\end{array}  
\end{eqnarray*}
This implies $v=0$. Therefore, $u=0$ and $w_1=w_2$. $\square$ \\


{\bf Corollary 2 (Conditions at interfaces).}
{\it Under the hypotheses of Theorem 1, if $h \in C^1([0,T];L^2(R))$,  
$u_1 \in H^1(R)$ and 
$u_2= {\rm div}(\chi \nabla u_0)/\rho + h(0,\mathbf x) \in L^2(R)$, 
the solution $u \in C^1([0,T];H^1(R)) \cap C^2([0,T];L^1(R))$ and additional estimates
\begin{eqnarray} 
\begin{array}{l} \displaystyle
 \|u_{tt}\|_{L^\infty(0,T;L^2(R))} \leq K(T, \rho_{\rm min}, \rho_{\rm \max}, E(u_1,u_2)), 
 \|h_t\|_{C([0,T];L^2(R))}),  \\[1ex]
 \| \nabla u_t \|_{L^\infty(0,T;L^2(R))} \leq
K(\chi_{\rm min}, T, \rho_{\rm min}, \rho_{\rm \max}, E(u_1,u_2), \|h_t\|_{C([0,T];L^2(R))}),  
 \\[1ex]
 \|u_{tt}\|_{L^2(0,T;L^2(\partial R \setminus \overline{\Sigma}))} \leq 
 K(\gamma_{\rm min}, \chi_{\rm min}, T, \rho_{\rm min}, \rho_{\rm \max}, E(u_1,u_2), 
 \|h_t\|_{C([0,T];L^2(R))}),  
\end{array} \label{continuous2}
\end{eqnarray}
hold. Normal derivatives of $u$ at the boundary are defined in $H^{-1/2}$ and
$u$ is a weak solution to (\ref{balance_escalar}).

Moreover, if $R = \cup_{\ell=1}^L R^\ell$,  $R^\ell$  being disjoint regions where 
$\rho$ and $\chi$ are constan, at each interface $\Gamma$ separating two 
adjacent regions
\begin{eqnarray}  \begin{array}{ll}
u^+ = u^-, & {\rm in } \, L^2(\Gamma), \\[1ex]
 \chi^-  \nabla u^- \cdot \mathbf n 
= \chi^+   \nabla u^+ \cdot \mathbf n, & {\rm  in } \, H^{-1/2}(\Gamma),
\end{array} \label{transmission}
\end{eqnarray}
where $^+$ and $^-$ denote limit values from each side
following the direction of unit normal vector $\mathbf n$.}

{\bf Proof.} Formally, differentiating  (\ref{balance_escalar}) and (\ref{balance_weak}) 
with respect to $t$, we have a similar variational equation for $v=u_t$ with right
hand side $h_t \in C([0,T];L^2(R))$ and  
$v(0, \mathbf x) = u_1 \in H^1(R)$ and 
$v_{t}(0, \mathbf x) = {\rm div}(\chi \nabla u_0)/\rho + h(0,\mathbf x)
= u_2 \in L^2(R)$. 
This problem admits a solution 
$v \in  C([0,T];H^1(R)) \cap C^1([0,T];L^2(R))$ satisfying the properties
stated in Theorem 1, which must be equal to $u_t$ in a distribution sense.
The enhanced regularity  implies that, for all t, 
$-{\rm div} (\chi \nabla u(t)) = \rho  h(t) -\rho  u_{tt}(t) \in L^2(R)$ with 
traces $u(t)|_{\partial R} \in L^2(\partial R)$.
Let us recall that $\mathbf q(t) \in L^2(R)$ and ${\rm div} (\mathbf q(t)) \in  L^2(R)$ 
imply $\mathbf q \cdot \mathbf n \in H^{-1/2}(\partial R)$,
see  \cite{bamberger}. Thus $\nabla u \cdot \mathbf n$ is defined at boundaries
in $H^{-1/2}$. Integrating by parts in (\ref{eq:v}) and using (\ref{eq:u})
we find 
\begin{eqnarray} \begin{array}{ll}
  \nabla v \cdot \mathbf n= - \gamma (v_t + \mu v)
\quad {\rm on} \; \partial R \setminus \overline{\Sigma}
& \implies  \nabla u \cdot \mathbf n = - \gamma u_t
\quad {\rm on} \; \partial R \setminus \overline{\Sigma}, \\[1ex]
 \nabla v \cdot \mathbf n= 0
\quad {\rm on} \; \Sigma
& \implies \nabla u \cdot \mathbf n = 0
\quad {\rm on} \; \Sigma.
\end{array} \label{eq:bcu}
\end{eqnarray}

Now, consider two adjacent domains  $R^1$ and  $R^2$ 
with common interface $\Gamma$.  Given $u \in H^1(R^1\cup R^2)$, and 
denoting by $u^+$, $u^-$ the  limit values of $u$ taken from 
$R^1$ and $R2$, respectively, we must have 
$u^+=u^-$ on $\Gamma$ in the sense of $L^2$ traces. Moreover, for 
any $w \in H^1_0(R^1\cup R^2)$ we have the identities
\begin{eqnarray*} \begin{array}{ll}
\sum_{\ell=1}^2 \int_{\Omega^\ell}    \chi^\ell  \nabla u  
\nabla w \, d \mathbf x  & 
 -  \sum_{\ell=1}^2 \int_{\Omega^\ell}  
 {\rm div}(\chi^\ell \nabla u) w \, d \mathbf x 
\\[1ex] 
& 
+ _{H^{-1/2}(\Gamma)}\langle [\chi^-  \nabla u^- \!-\!
\chi^+  \nabla u^+] \cdot \mathbf n, w \rangle_{H^{1/2}(\Gamma)},
\\[1ex] 
\int_{\Omega^1\cup \Omega^2}    \chi \,
 \nabla u   \nabla w \, d \mathbf x  & 
 = -\int_{\Omega^1\cup \Omega^2}  
{\rm div}(\chi \nabla u) w \, d \mathbf x.
\end{array} 
\end{eqnarray*}
Therefore,  the transmission relations hold at any discontinuity 
interface $\Gamma$.  $\square$ \\


{\bf Theorem 3 (Regularity).} {\it Under the hypotheses of Theorem 1 and 2, let us further
assume that
\begin{itemize}
 \item $u_0=u_1=0$,
 \item $h(t, \mathbf x) = f(t) g(\mathbf x) \in C^2([0,T];H^1(R))$, where $f\in C^2([0,T])$
and $g(\mathbf x) \in H^1(R)$, with ${\partial g \over \partial \mathbf n} = 0$ on $\Sigma$,
support contained in an upper layer $R^1$ with upper boundary $\Sigma$, a
nd vanishing at a positive distance of ${\partial R^1 \setminus \Sigma}$.
\end{itemize}
Then, the solution of (\ref{balance_escalar})  satisfies 
$u \in C^2([0,T]; H^1(R)) \cap C^3([0,T]; L^2(R)) $ to problem (\ref{balance_escalar})  
with
$u_t, u_{tt}, u_{ttt} \in L^2(0,T;L^2(\partial R \setminus \overline{\Sigma}))$. 
The following additional stability estimates hold 
\begin{eqnarray} \begin{array}{l} \displaystyle
\|u_{ttt}\|_{L^\infty(0,T;L^2(R))} \leq K(T, \rho_{\rm min}, \rho_{\rm \max}, E(u_2,u_3), 
\|h_{tt}\|_{C([0,T];L^2(R))}),  \\[1ex]
 \| \nabla u_{tt} \|_{L^\infty(0,T;L^2(R))} \leq
K(\chi_{\rm min}, T, \rho_{\rm min}, \rho_{\rm \max}, E(u_2,u_3), \|h_{tt}\|_{C([0,T];L^2(R))}),  
\\[1ex]
 \|u_{ttt}\|_{L^2(0,T;L^2(\partial R \setminus \overline{\Sigma}))} \leq 
K(\gamma_{\rm min}, \chi_{\rm min}, T, \rho_{\rm min}, \rho_{\rm \max}, E(u_2,u_3),
\|h_{tt}\|_{C([0,T];L^2(R))}), 
\end{array} \label{continuous2} \end{eqnarray}
where $u_2=f(0)g$, $u_3=f'(0)g$ and $K(\cdot)$ denote different positive
 constants depending continuously on the specified arguments. If $g\in H^2(R)$ 
 and $f\in C^3([0,T])$,
analogous regularity and estimates hold for $u_{tttt}$.

Moreover, $u$ has $H^2$ regularity in the upper layer near $\Sigma$. In dimension
$n=2$,  $u$ is continuous up to $\Sigma$  and its values at the receptor
points  located at $\Sigma$   are defined, at least during a certain time.
}

{\bf Proof.} 
Differentiating (\ref{balance_weak}) twice with respect to $t$, we find for $v=u_{tt}$ 
a similar problem, with right hand side $h_{tt} \in C([0,T];L^2(R))$ and initial data 
$v(0, \mathbf x) = f(0) g(\mathbf x) \in H^1(R)$ and $v_{t}(0, \mathbf x) = f'(0)
g(\mathbf x) \in L^2(R)$. This yields   (\ref{continuous2}).

Moreover, for each fixed $t>0$,  we have
\begin{eqnarray*}\begin{array}{ll}
{\rm div}(\chi(\mathbf x) \nabla u(t)) = \rho(\mathbf x)  h(t,\mathbf x)
- \rho(\mathbf x) u_{tt}(t) \in L^2(R), & \mathbf x \in R,\\[1ex]
\nabla u(t) \cdot \mathbf n = 0,  & \mathbf x \in \Sigma, \\ [1ex]
\nabla u(t) \cdot \mathbf n = - \gamma(\mathbf x) u_{t} \in L^2(\partial R 
\setminus \overline{\Sigma}), & \mathbf x \in \partial R \setminus \overline{\Sigma}.
\end{array}\end{eqnarray*}
Let us consider a smooth function $\eta(x_1,x_2)$ that decreases from $1$ when 
$x_2 \in [0,-\delta]$ to $0$ at $x_2 =- 2 \delta$, $\eta>0$ small enough, with 
support of $\eta$ to be contained in the upper layer $R^1$ around the receivers.
Then,  $w=  u \eta$ has support contained in the first layer and satisfies
\begin{eqnarray*}\begin{array}{ll}
\chi \Delta w = \eta [\rho  h(t) - \rho u_{tt}(t)] 
+ 2 \chi \nabla u \nabla \eta + \chi u \Delta \eta = r, & \mathbf x \in R, \\[1ex]
\nabla w \cdot \mathbf n = 0,  & \mathbf x \in \Sigma, \\ [1ex]
\nabla w \cdot \mathbf n = - \eta \gamma u_{t} + u \nabla \eta \cdot
\mathbf n & \mathbf x \in \partial R \setminus \overline{\Sigma}, \\[1ex]
\end{array}\end{eqnarray*}
and also, $ w = 0, \mathbf x \in \Sigma_b, $ where $\Sigma_b$ is the bottom wall. 
For any $g \in L^2(\partial R)$, there exists a function $w_g \in H^2(R)$ such that
$\nabla w_g \cdot \mathbf n = g$. Setting $w = \tilde w + w_g$, we find
\begin{eqnarray*}\begin{array}{l}
\chi  \Delta \tilde w(t)  =  r  -  \chi \Delta w_g(t) \in L^2(R), 
\\ [1ex]
\nabla \tilde w(t) \cdot \mathbf n = 0,  \quad {\rm on }  \, \partial R.
\end{array}\end{eqnarray*}
By the regularity results in \cite{grisvard}, Ch. 3, the solution $\tilde w \in H^2(R)$. 
Therefore, $u$ has $H^2$ regularity in a neighborhood of the receivers for a certain 
time. By Sobolev's injections \cite{grisvard}, $u$ is continuous up to the border in 
that region in dimension $n=2$. In particular, it is defined at receptor points. $\square$ 


Let us discuss now how to approximate numerically the solution of the
forward problem.


\section{Appendix C: Numerical approximation of the truncated forward problem}
\label{sec:truncated_na}

Applications to inverse problems require solving large amounts of wave problems, 
therefore, its is desirable to keep the computational cost as low as possible.
We discuss here the approximations employed together with their stability and
convergence properties.


\subsection{Space-time discretization}
\label{sec:discretization}

We resort to finite element dicretizations in space and finite differences
in time. Given a spatial mesh and an associated triangulation ${\cal T}_{\delta x}$,
with maximum time step $\delta x$, we build a FEM basis $\{\psi_1, \ldots, \psi_{M}\}$, 
$M=M(\delta x)$, of $P^1$ elements at least. 
Let $V^{M}$ be the space spanned by them.

The discretized problem  becomes: Find
$u^M(t,\mathbf x) = \sum_{m=1}^M a_m(t) \psi_m(\mathbf x)$ 
such that
\begin{eqnarray} 
\begin{array}{l} 
{d^2\over dt^2} \int_{R}  \rho u^M(t)  w  \, d \mathbf x +
\int_{R}  \chi \nabla  u^M(t)  \nabla  w  \, d \mathbf x   
+ {d \over dt} \int_{\partial R \setminus \overline{\Sigma}} 
\chi \gamma \, u^M(t)  w \, d S_{\mathbf x}  \\
 \hfill =   f(t)  \int_{R} \rho g  w \,  d \mathbf x, 
\\[2ex] 
u^M(0) =  0, \quad u_t^M(0) =  0, 
\end{array}  \label{eq:femM}
\end{eqnarray}
for all $w \in V^M$  and $t \in [0,T]$.

{\bf Lemma 4 (Existence of an approximant). }
{\it For $f \in C^2([0,T])$ there is a unique 
function  $u^{M} \in C^2([0,T]; V^{M})$ satisfying (\ref{eq:femM}). 
Moreover, if $f \in C^k([0,T]$, then $u^{M} \in  C^{k+2}([0,T]; V^{M})$,
$k\geq 1$.} \\
{\bf Proof.}
Taking $w=\psi_j$, $j=1,\ldots,M$, (\ref{eq:femM}) is equivalent
to a system of ordinary differential equations
\begin{eqnarray}
\begin{array}{r} 
\sum_{m=1}^M a_m''(t) \int_{R} \rho  \psi_m  \psi_j \, d \mathbf x  
+ \sum_{m=1}^M a_m(t) \int_{R}  \chi
\nabla \psi_m \nabla \psi_j \, d \mathbf x 
\\[2ex] 
+ \sum_{m=1}^M a_m'(t) \int_{\partial R \setminus \overline{\Sigma}} 
\chi \gamma \, \psi_m  \psi_j \, d S_{\mathbf x} 
=   f(t) \int_{R} \rho g(\mathbf x)  \psi_j\, d \mathbf x, 
\\[2ex] 
a_m(0)= a_m'(0) = 0, \quad m=1, \ldots, M.
\end{array}  \label{discrete_system}
\end{eqnarray}
In vector form
$\mathbf B \mathbf v'(t) = - \mathbf C  \mathbf a(t)  - \mathbf E
\mathbf v(t) + \mathbf f(t) \mathbf h, $
$\mathbf a'(t)  =  \mathbf bv(t),  $
$\mathbf v(0) = \mathbf a(0) = 0,$ 
$t \in [0, T]$.
Notice that $\mathbf M$ is symmetric and positive definite. This linear
system of differential equations  has a unique
solution in $[0,T]$, see \cite{coddington}, whose regularity increases
with the regularity of $f(t)$. $\square$ \\

One can also work with more refined variational formulations
set in the different domains and connected through an additional
bilinear form representing the transmission conditions \cite{bamberger,
chabassier}. 

Let us consider now the time discretization. We discretize the time derivatives 
in (\ref{eq:femM}) using centered differences for $u_{tt}^{M}$ 
\begin{eqnarray*} \begin{array}{l}
u_{tt}^{M}(\mathbf x,t) \sim {u^{M}(\mathbf x,t+\delta t) - 2u^{M}(\mathbf x,t)
+ u^{M}(\mathbf x,t-\delta t) \over \delta t^2 } + O(\delta t^2), 
\end{array} \end{eqnarray*}
and backward or centered differences for $u_{t}^{M}$
\begin{eqnarray} 
\begin{array}{l}
u_{t}^{M}(\mathbf x,t)  \sim {u^{M}(\mathbf x,t) - u^{M}(\mathbf x,t-\delta t) 
\over \delta t} + O(\delta t), 
\end{array}  \label{1st} \\
\begin{array}{l}
u_{t}^{M}(\mathbf x,t)  \sim {u^{M}(\mathbf x,t+\delta t) - u^{M}(\mathbf x,t-\delta t) 
\over 2\delta t} + O(\delta t^2).  
\end{array} \label{2nd}
 \end{eqnarray}
The first choice yields the scheme
\begin{eqnarray}
\begin{array}{r} 
\int_{R} \rho \,   u^{M}(t+\delta t)  w\, d \mathbf x   
=  \int_{R}  \rho \,   [2u^{M}(t)-u^{M}(t-\delta t)]  w \, d \mathbf x  
\\[1.5ex]  
- \delta t^2 \int_{R}  \chi  \, \nabla u^{M}(t) \, \nabla w\, d \mathbf x  
+ \delta t^2  \int_{R}   \rho  f(t) w \, d \mathbf x 
\\[1.5ex]  
- \delta t \int_{\partial R - \Sigma} \rho \chi \gamma
\, [u^{M}(t)-u^{M}(t-\delta t)] w \, d S_{\mathbf x}.
\end{array} \label{time_artificial}
\end{eqnarray}
On a temporal grid $t_n=n\, \delta t$, $n=0,\ldots,N$, $n\delta t = T$ we approximate 
$u^M(t_n) = \sum_{m=1}^M a_m(t_n) \phi_m$by $\sum_{m=1}^M a_m^n \phi_m$.
The coefficients $a_m(t_n)$, $m=1,\ldots M$, are approximated by the
solution $a_m^n$ of the recurrence:
\begin{eqnarray} \begin{array}{r} \displaystyle
\sum_{m=1}^{M} B_{j,m} a_m^{n+1} =  
\sum_{m=1}^{M}  B_{j,m} (2 a_m^n  - a_m^{n-1}) 
-  \delta t^2  \sum_{m=1}^M C_{j,m} a_m^n
\\ [1ex] \displaystyle 
-  \, \delta t \sum_{m=1}^M E_{j,m}   (a_m^n- a_m^{n-1})  
+ \delta t^2  f(t_n) h_j, \quad j=1, \ldots, M,
\end{array} \label{discretization} \end{eqnarray}
for $m \geq 1$, with $\mathbf B, \mathbf C, \mathbf E, \mathbf h$ defined by
\begin{eqnarray} \begin{array}{ll}
B_{j,m} = \int_{R} \rho  \psi_m  \psi_j \, d \mathbf x,  &
C_{j,m} = \int_{R}  \chi \nabla \psi_m \nabla \psi_j \, d \mathbf x, \quad 
\\[2ex] \displaystyle
E_{j,m} = \int_{\partial R \setminus \Sigma}  \chi \gamma \psi_j \psi_m 
\, d S_{\mathbf x},  &
h_j = \int_{R} \rho g(\mathbf x)  \psi_j\, d \mathbf x,
\end{array} \label{matrices_artificial}
\end{eqnarray}
for $j,m=1,\ldots, M$.
Initially, $a_m^0=a_m(t_0)=0$ and $a_m^1=a_m(t_0)+\delta t \,
a_m'(t_0)=0$ for $m=1,\ldots,M$, as dictated by the initial conditions. 

{\bf Lemma 5 (Existence of an approximant). }
{\it Scheme  (\ref{discretization}) admits a unique solution 
$a_m^n$, $n=0,\ldots,N$, $m=1,\ldots,M$. } \\
{\bf Proof.} Given two time levels $a_m^{n-1}$ and $a_m^{n}$,
$n=0,\ldots,N$, the level $n+1$ follows directly from relation
 (\ref{discretization}). $\square$ 
 
The matrices and vectors involved in system (\ref{discretization}) are
known and fixed once the spatial mesh and the associated finite element 
basis are constructed, provided $\rho$, $\chi$, $\gamma$ remain unchanged. 
For a given inclusion, we use the same matrices and vectors in all the temporal
steps. If we vary the shape and the material parameters of inclusions, we 
need to recalculate them. Depending on whether we keep the mesh fixed
or update it, we need to remesh and calculate a new function basis too.


\subsection{Convergence results}
\label{sec:convergence}


Consider a regular family of triangulations ${\cal T}_{\delta x}$ of $R$ with 
maximum element diameter $\delta x \rightarrow 0$. Let us define associate
$P^1$ finite element spaces $V^M\subset H^1(R)$ of dimension 
$M=M(\delta x) \rightarrow \infty$. We introduce the elliptic projection operator 
$\Pi^M: H^1(R) \rightarrow V^M$ that associates to each $v \in H^1(R)$
the solution $\Pi^M v \in V^M$ of the elliptic problem
\begin{eqnarray} \begin{array}{l}
\displaystyle
a_\mu(v,v^M) = a_\mu(\Pi^M v, v^M), \quad \forall v^M \in V^M, \\
[1ex] \displaystyle 
a_\mu(v,w) = \int_R \chi \nabla v \nabla w d \mathbf x + 
\int_R \rho \mu^2  \, v w d \mathbf x + \int_{\partial R \setminus \overline{\Sigma}}
\chi \gamma \mu \, v w dS_{\mathbf x},
\end{array} \label{projection}
\end{eqnarray}
associated to (\ref{eq:varwM_lam}). 

{\bf Theorem 6 (Convergence of the FEM discretization).} 
{\it Consider the solution $u \in C^2([0,T];H^1(R))$ of (\ref{balance_escalar})
constructed under the hypotheses of Theorems 1-3 and 
$u^M \in C^2([0,T];V^M) \subset H^1(R)$ the solutions generated by the 
FEM scheme (\ref{eq:femM}).
Assume that the family of regular triangulations considered satisfies the 
approximation property
\begin{eqnarray}
{\rm lim}_{M \rightarrow \infty} {\rm inf}_{v^M \in V^M} \| v - v^M \|_{H^1(R)}=0.
\label{approximation}
\end{eqnarray}
Then, the sequence $(u^M(t),u^M_t(t))$ converges to $(u(t),u_t(t))$ in 
$H^1(R) \times L^2(R)$ for all $t \in [0,T]$.}

{\bf Proof.}  We make the change of variables $u= e^{\mu t} v$ and
$u^M = e^{\mu t} v^M$ in the corresponding variational equations.
Subtracting (\ref{projection}), we find \cite{raviart}
\begin{eqnarray*} \begin{array}{l}
{d^2\over dt^2} \int_{R}  \rho (v^M-\Pi^M v)(t)  v^M  \, d \mathbf x 
+ {d \over dt} \int_{R}  2 \rho  \mu  (v^M-\Pi^M v)(t)  v^M  \, d \mathbf x
\\[1ex]
+ {d \over dt} \int_{\partial R \setminus \overline{\Sigma}} 
\chi \gamma \, (v^M-\Pi^M v)(t)  v^M \, d S_{\mathbf x} +
a((v^M-\Pi^M v)(t), v^M)  
\\[1ex]
=  \int_{R}  \rho (I-\Pi^M) v_{tt}(t)  v^M  \, d \mathbf x 
+ \int_{R}  2 \rho  \mu  (I-\Pi^M) v_t(t)  v^M  \, d \mathbf x
\\[1ex]
+ \int_{\partial R \setminus \overline{\Sigma}} 
\chi \gamma   (I-\Pi^M)  v_t(t)  v^M \, d S_{\mathbf x} 
\end{array} \label{eq_projection}
\end{eqnarray*}
since $v \in C^2([0,T];H^1(R))$, with zero initial data. The solution
$v^M-\Pi^M v$ of this problem satisfies an energy inequality
analogous to (\ref{eq:ienergyv}). Applying Young's inequality to
the right hand side, and taking into account the zero initial data,
we find
\begin{eqnarray*} \begin{array}{l} 
E_\mu((v^M-\Pi^M v)(t), (v^M-\Pi^M v)_t(t)) \leq 
C(\rho,\mu,\chi) \int_0^t  [ \| (I-\Pi^M) v_{tt}(t) \|^2_{L^2(R)} 
\\ [2ex] 
+ \| (I-\Pi^M) v_{t}(s) \|^2_{L^2(R)} +
\| (I-\Pi^M) v_{t}(s) \|^2_{L^2(\partial R \setminus \overline{\Sigma})} ]ds 
= I^M(t)
\end{array} \label{energy_projection}
\end{eqnarray*}
with $E_\mu$ defined in (\ref{en_def_mu}). Then,
$E_\mu((v-\Pi^M v)(t), (v-v^M)_t(t))$ is bounded from above in terms of
$E_\mu((v-\Pi^M v)(t), (v-v^M)_t(t))$ and $I^M(t)$.
Condition (\ref{approximation}) implies that
$E_\mu((v-\Pi^M v)(t), (v-v^M)_t(t)) \rightarrow 0$ as $M \rightarrow \infty$.
Undoing the change, the same is true for $u$, therefore, the FEM
approximation converges to the solution. 
Notice that we work in a polygonal domain whose external boundary
is fixed. $\square$ 

{\bf Remark 1.} For functions $v \in H^2(\Omega)$, estimates
of the form $\| v - v^M \|_{H^1(R)} \leq C \delta x \| v \|_{H^2(R)}$
hold for regular triangulations \cite{raviart} and $P^1$ elements, which
ensures condition (\ref{approximation}) and $o(\delta x)$ convergence. 
Solutions of (\ref{balance_escalar})  with piecewise constant coefficients 
in subdomains $R^\ell$, $\ell )a, \ldots,L$ can reach at best
$H^2(R^\ell)$ regularity in each subdomain. No global $H^2(R)$ regularity
can be achieved. In general, we can expect convergence when triangulations 
are made following a domain decomposition approach, that is, they are 
entirely contained in each region $R^\ell$, sharing vertices at the interfaces. 

{\bf Remark 2.} An additional issue regarding convergence of FEM approximations
stems from the fact that we must approximate numerically the integrals 
appearing in (\ref{eq:femM}). Convergence of these approximations
is easier to control when the triangulations are contained in each subdomain
$R^\ell$ sharing nodes at interfaces, avoiding triangles partially contained
in different subdomains and mesh nodes moving from one subdomain to
another as the triangulations are refined. \\

{\bf Lemma 7 (Convergence of the time discretization).}
{\it Consider the solution $u^M=\sum_{m=1}^M a_m(t) \phi_m \in C^2([0,T])$  
of (\ref{eq:femM}) constructed in Lemma 4. The sequence $u_m^n$, 
$n=0, \ldots, N$ generated by scheme (\ref{discretization}) converges
in the sense that
${\rm max_{n=0,\ldots,N}} \| a_m(t_n) - u_m^n \| \rightarrow 0$
as $N \rightarrow \infty$.} 

{\bf Proof.}
The approximation has truncation error $O(\delta t^2+ \delta t)$.
Convergence requires that matrices
$2\mathbf I  -  \delta t^2 \mathbf B^{-1} \mathbf C  - \delta t \mathbf B^{-1} \mathbf E$
and $\mathbf I- \delta t \mathbf B^{-1} \mathbf E$ have spectral radius smaller than
$1$. These matrices depend on $\delta x$, which determines their size $M$, thus
$\delta t/\delta x$ must remain small enough \cite{IKeller,raviart}. \\

{\bf Corollary 8 (Convergence of the full discretization).} 
{\it Let $u \in C^2([0,T];H^1(R))$ be the solution of (\ref{balance_escalar})
constructed under the hypotheses of Theorems 1-3 and 
$u^M_n \in V^M$, $n=0,\ldots, N$, the sequences generated by the 
scheme (\ref{discretization}). Then the error 
${\rm max_{n=0,\ldots,N}} \| u^M(t_n) - u^M_n \|_{H^1}(R)$
tends to zero as $N,M \rightarrow \infty$ provided
\begin{itemize}
\item condition (\ref{approximation}) holds in the regular triangulation
uniformly for $v$ belonging to a bounded $H^1(R)$ set,
\item the ratio ${\delta t \over \delta x}$ of the time step $\delta t$ to the
spatial diameter $\delta x$ of the triangulation is small enough.
\end{itemize}}
{\bf Proof.}
Consequence of Theorem 6 and Lemma 7. $\square$\\

{\bf Remark 3.} When $R=\cup_{\ell=1}^L R^\ell$, a Courant-Friedrichs-Lewy 
condition of the  form $\delta t \leq {\rm Min}_{\ell=1,\ldots, L}\{1/v_{p,\ell}\} 
\delta x/2 $ in terms of the known layer waves speeds \cite{john} preserves 
stability in the simulations performed here.
In the simulations shown here, we have used the second
order discretization (\ref{2nd}) in the boundary terms, increasing the
global approximation order, and still preserving stability while slightly
modifying the scheme. 


\end{document}